\documentclass[12pt]{amsart}
\usepackage{amsmath}
\usepackage{amssymb}
\usepackage{tabularx}
\usepackage{enumerate}
\usepackage{graphicx}
\usepackage{texdraw}
\usepackage{exscale}
\usepackage{relsize}
\usepackage{booktabs}
\usepackage{multirow}
\usepackage{mathrsfs}
\usepackage{latexsym}
\usepackage{amsfonts,amssymb,amsmath}
\usepackage{epsfig}
\usepackage{color}
\makeatletter
\@addtoreset{equation}{section}

\makeatother
\newtheorem{thm}{Theorem}[section]
\newtheorem{lem}[thm]{Lemma}

\newtheorem{prop}[thm]{Proposition}
\newtheorem{remark}[thm]{Remark}

\newcommand{\R}{\mathbb{R}}

\begin{document}
\title[A reaction-diffusion model with Allee effect and free boundary]
{Long-time behavior of a reaction-diffusion model with strong Allee effect and free boundary:\ effect of a protection zone$^\S$}

 \thanks{$\S$ N. Sun was  partially supported by NSF of China (No. 11801330),
 the Higher Educational Science and Technology Program of Shandong Province (No. J18KA226), and  C. Lei was partially supported by
 NSF of China (No. 11801232), the Priority Academic Program Development of Jiangsu Higher Education Institution, the Natural Science
 Foundation of the Jiangsu Province (No. BK20180999), the Foundation of Jiangsu Normal University (No. 17XLR008).}
\author[N. Sun, C. Lei]{Ningkui Sun$^\dag$ and Chengxia Lei$^\ddag$}
\thanks{$\dag$ School of Mathematics and Statistics, Shandong Normal University, Jinan, 250014, Shandong Province, China.}
\thanks{$\ddag$ School of Mathematics and Statistics, Jiangsu Normal University, Xuzhou, 221116, Jiangsu Province, China.}
\thanks{E-mails: {{\sf sunnk@sdnu.edu.cn} (Sun), \sf leichengxia001@163.com} (Lei)}
%\thanks{$*$ Corresponding author.}
\date{}

\begin{abstract}
This paper concerns the effect of the (separated/connected) protection zone for the evolution of an endangered species on the
reaction-diffusion equation with strong Allee effect and free boundary. We give a description
of the long-time dynamical behavior of the problem of two types protection zones with the same length. Furthermore, the asymptotic profiles of solutions
and the asymptotic spreading speed are estimated when spreading happens. Our results, together with those in previous papers \cite{DPS, DuLou} on two other closely related
models, show that the protection zone and the free boundary play an important role in the evolution of the endangered species.

\end{abstract}

\maketitle

\section{Introduction}\label{sec:pro}
More than 99 percent of all species amounting to over five billion species, that ever lived on Earth are estimated to be extinct;
most species that become extinct are never scientifically documented \cite{Ne}. Some scientists estimate that up to half of presently
existing plant and animal species may become extinct by 2100 \cite{Wi}. There are more and more environmental groups
and governments paying attention to the endangered species. A variety of protection planning projects have been taken to protect the
endangered species and their habitats. One of the most effective strategies of preventing the endangered species from reducing to
lower densities or small number is constructing the protection zone. The effect of protection zones on endangered species has attracted
much attention by the biologists and the mathematicians. Recently, some mathematical models describing protection zones have been
proposed and analyzed; see, for instance, \cite{CSW,DPS,DuLiang,DPW,DuS2,HZ,LWL,Oe} and the references therein.

In particular, the authors in \cite{DPS} introduced a reaction-diffusion model with strong Allee effect and a protection
zone to examine the role of the protection zone on the dynamical behavior of the species evolution. They considered the following problem
\begin{equation}\label{20201031}
\left\{
\begin{array}{ll}
u_t = u_{xx} + f(u), &  t>0, -L<x<L,\\
u_t = u_{xx} + g(u), &  t>0, x\in \mathbb{R} /[-L,L],\\
u(t, L-0)= u(t, L+0), &t>0,\\
u(t, -L-0)= u(t, -L+0)& t>0,\\
u_x(t, L-0)= u_x(t, L+0),& t>0,\\
u_x(t, -L-0)= u_x(t, -L+0),&t>0,\\
u(0,x)=u_0(x)\geq 0, & x\in \mathbb{R},
\end{array}
\right.
\end{equation}
where $L>0$, and the protection zone is $[-L, L]$; $u(t,x)$ stands for the population density
of the species under consideration at time $t$ and location $x$; the initial function $u_0(x)$ is nonnegative and compactly supported; for any given $t>0$,
$u(t, L-0)$ and $u_x(t, L-0)$ represent, respectively, the left limit value and the left
derivative of $u$ with respect to $x$ at $x=L$, and $u(t, L+0)$ and $u_x(t, L+0)$ are the right limit value and the right derivative of $u$ with respect to $x$ at
$x=L$, respectively.

In the protection zone, the growth of the species is governed by a monostable nonlinearity
$f(u)$ which satisfies
 \begin{equation}\label{mono}
 f(0)=f(1)=0<f'(0), \ \ f'(1)<0,\ \ (1-u)f(u) >0,\ \ \forall u>0, \ u\neq 1.
 \end{equation}
The nonlinearity $g(u)$ is used to describe the evolution species which obeys the strong Allee effect \cite{A} out of the protection zone.
To include the strong Allee effect, a typical reaction function is the so-called ``bistable'' nonlinear terms; see, for example,
\cite{CBG,KLH,LKa,ML,WSW} and the references therein. The bistable nonlinear term $g$ is {\it globally Lipschitz}
and satisfies
 \begin{equation}\label{bi}
 g(0)=g(\theta)= g(1)=0, \quad g(u) \left\{
 \begin{array}{l}
 <0 \ \ \mbox{in } (0,\theta),\\
 >0\ \  \mbox{in } (\theta, 1),\\
 < 0\ \ \mbox{in } (1,\infty),
 \end{array} \right.
 \end{equation}
 for some $\theta\in (0,1)$,  $g'(0)<0$, $g'(1)<0$ and
 \begin{equation}\label{unbalance}
 \int_0^1 g(s) ds >0.
 \end{equation}

For simplify, they  investigated the following symmetrical form:
 \begin{equation}\label{pp1}
 \left\{
 \begin{array}{ll}
 u_t = u_{xx} + f(u), &  t>0, 0<x<L,\\
 u_t = u_{xx} + g(u), &  t>0, x>L,\\
 u_x(t,0)=0,&t>0\\
 u(t, L-0)= u(t, L+0),& t>0,\\
 u_x(t, L-0)= u_x(t, L+0),& t>0,\\
 u(0,x)=u_0(x)\geq 0, & x\geq 0.
 \end{array}
 \right.
 \end{equation}
It was shown in \cite{DPS} that there are two critical values $0<L_*\leq L^*$ which affect the dynamics of the solutions
significantly. More precisely, in the small protection zone case($L<L_*$), there is a vanishing-transition-spreading
trichotomy result;
in the medium-sized protection zone case($L_*<L<L^*$), there is a transition-spreading dichotomy result;
only {spreading} happens in the large protection zone case($L>L^*$). They  compared two types of protection zone with the same length:
a connected one and a separate one, and their results revealed that the former is better for species
spreading than the latter.

When $L=0$, that is to say, there is no protection zone in the environment,
 some special cases of \eqref{20201031} were studied by many authors. Among them,
Du and Matano \cite{DM} obtained a rather comprehensive analysis of the dynamical behavior of solutions
by introducing a parameter in the initial value. And any other relevant works can be found, for instance,
in \cite{AW2, Zla} and the references therein. In the case of $L=\infty$, the authors in \cite{AW2} also studied problem  \eqref{20201031},
and obtained that the asymptotic behavior of solutions and the existence of the traveling wave solutions.
For the special logistic nonlinearity $f(u)=u(1-u)$, it was also focused by \cite{F,KPP} to describe the spreading of
an advantageous genetic trait in a population.

In addition, for this logistic nonlinearity $u(1-u)$, a free boundary problem
\begin{equation*}\label{freeb}
\left\{
\begin{array}{ll}
u_t =u_{xx}+u(1-u), & t>0,\; 0< x<h(t),\\
u_x(t,0)=0=u(t,h(t)), & t>0,\\
h'(t)=-\mu  u_x  (t, h(t)), & t>0,\\
h(0)=h_0,\;u(0,x)=u_0(x)\geq 0, &  0\leq x \leq h(0),
\end{array}
\right.
\end{equation*}
was proposed by Du and Lin \cite{DuLin} to describe the spreading of a new or invasive species in which the free boundary $x = h(t)$
represents the spreading front of the population whose density is represented by $u(t, x)$, the positive coefficient $\mu$ measures the
ability for spreading into new habitat of the species (for more background, see \cite{BDK, L}). Later on, Du and Lou \cite{DuLou} studied the following two free boundary problems
\begin{equation}\label{20201031-1}
\left\{
\begin{array}{ll}
u_t =u_{xx}+f(u), & t>0,\; 0< x<h(t),\\
u_x(t,0)=0=u(t,h(t)), & t>0,\\
h'(t)=-\mu  u_x  (t, h(t)), & t>0,\\
h(0)=h_0,\;u(0,x)=u_0(x)\geq 0, &  0\leq x \leq h_0,
\end{array}
\right.
\end{equation}
and
\begin{equation}\label{20201031-2}
\left\{
\begin{array}{ll}
u_t =u_{xx}+g(u), & t>0,\; 0< x<h(t),\\
u_x(t,0)=0=u(t,h(t)), & t>0,\\
h'(t)=-\mu  u_x  (t, h(t)), & t>0,\\
h(0)=h_0,\;u(0,x)=u_0(x)\geq 0, &  0\leq x \leq h_0.
\end{array}
\right.
\end{equation}
For further related work on free
boundary problems, we refer to \cite{DuGuo,DMZ2, S, SH, SF, W1} and the references therein.

By comparing the results of  problem \eqref{pp1} with $L=\infty$ and the corresponding free boundary problem of \eqref{20201031-1},
the striking difference between those two problems was that there is a ``hair-trigger'' phenomenon (namely $u(t,x)\to 1$
as $t\to\infty$) for the former problem \eqref{pp1} with any nontrivial initial value, regardless of its initial size and
supporting area any nontrivial initial value, but the dynamics of the free boundary problem exhibits a
spreading-vanishing dichotomy, please see \cite{AW2, DuLin} for more details.
The authors in \cite{DuLou} pointed out that problem \eqref{pp1}
with $L=\infty$ is the limiting problem
of the corresponding free boundary problem of \eqref{20201031-1}. Moreover, the asymptotic spreading speed of the free boundary problem is increasing to
the asymptotic spreading speed of the problem \eqref{pp1}
with $L=\infty$ as $\mu\rightarrow\infty$.
The existence of free boundary makes the dynamics of the solutions of \eqref{20201031-1} much more
complicated and meaningful than that of problem \eqref{pp1} with $L=\infty$. Moreover, we see that the presence of free boundary makes the
species more difficult to survive, please see the last discussion section for more details.

In the current paper, we consider the evolution of the species is governed
by the monostable nonlinear terms in the (separated/connected) protection zone, but the species growth obeys strong Allee effect out of the protection zone
and the expanding front is determined by a free boundary. Our model with the connected protection zone ($x\in[0,L]$) is described by the following form:
\begin{equation}\label{p}
\left\{
\begin{array}{ll}
u_t = u_{xx} + f(u), &  t>0, 0<x<L,\\
u_t = u_{xx} + g(u), &  t>0, L<x<h(t),\\
u_x(t, 0)=0= u(t,h(t)), & t>0,\\
h'(t)=-\mu u_x(t,h(t)),& t>0,\\
u(t, L-0)= u(t, L+0), & t>0,\\
u_x(t, L-0)= u_x(t, L+0), & t>0,\\
h(0)=h_0>L,\ u(0,x)=u_0(x)\geq 0, & x\in [0,h_0],
\end{array}
\right.
\tag{$P$}
\end{equation}
and the separated one ($x\in[L_1,L_2]$ with $L_2>L_1>0$ and $L_2-L_1=L$) reads as
\begin{equation}\label{q}
\left\{
\begin{array}{ll}
u_t = u_{xx} + f(u), &  t>0,  x \in(L_1, L_2),\\
u_t = u_{xx} + g(u), &  t>0, x\in(0, L_1)\cup(L_2,h(t)),\\
u(t, L_i-0)= u(t, L_i+0), & t>0, i=1, 2, \\
u_x(t,  L_i-0)= u_x(t, L_i+0), & t>0, i=1, 2, \\
u_x(t,0)=u(t,h(t))=0, & t>0,\\
h'(t)=-\mu u_x(t,h(t)), & t>0,\\
h(0)=h_0>L_2,\ u(0,x)=u_0(x)\geq 0, & x\in [0,h_0],
\end{array}
\right.\tag{$Q$}
\end{equation}
where the length of the protection zone is $L_2-L_1=L$, $x = h(t)$ stands for the spreading front of the population
whose density is represented by $u(t, x)$, $f$ is a monostable nonlinearity satisfying \eqref{mono},
and $g$ is a bistable nonlinearity satisfying \eqref{bi}, \eqref{unbalance}. Furthermore, we
assume that

\smallskip

{\bf{(H)}}\ \ \ \ \ \ The functions $f,\,g$ are {\it globally Lipschitz} and
$g(u)<f(u)\ \mbox{ for all } 0 < u < 1$.

\smallskip

The initial function $u_0$ belongs to  $ \mathscr {X}(h_0)$, where
\begin{equation*}\label{def:X}
\begin{array}{ll}
\mathscr {X}(h_0):= \big\{ \phi \in C^2 ([0,h_0]):
\phi'(0)=\phi (h_0)=0,\;
\phi(x)>0 \ \mbox{in } (0,h_0)\big\}.
\end{array}
\end{equation*}

In this paper, we have two main goals.  One of the objective is to investigate the long-time dynamical
behavior of problems \eqref{p} and \eqref{q} and the effect of the free boundary conditions and the structure of the protection
zone on the evolution of the endangered species.  Theorems \ref{thm:dybe} and \ref{thm:dybe1} in this paper tell us that only if the protection
zone is suitably long, the endangered species will survive in the entire space regardless of its initial data. Though similar conclusions
were also found of problem \eqref{pp1} in \cite{DPS}, the dynamics of solutions of problems \eqref{p} and \eqref{q} are much more complicated
and meaningful than those in \cite{DPS}. The other purpose  is to compare our results on the problem \eqref{p} with those on problems \eqref{pp1} and \eqref{20201031-2},
in order to comprehend the influence of the factors such as the free boundary condition and the protection zone on the persistence of the endangered species.
A more detailed description of the analytical results and the comparisons between three related problems
mentioned above can be found in the last section.

\smallskip

Now, we are going to state our main results of problems \eqref{p} and \eqref{q}.
Firstly, for any given $h_0>L$ and $u_0 \in \mathscr {X}(h_0)$, it follows from Section 2 that
\eqref{p} admits a unique time-global solution $(u(t,x), h(t))$ satisfying for any $\alpha\in(0,1)$,
\[
u\in C^{1,2}([0,\infty)\times([0,h(t)]/\{L\}))\cap C^{\frac{1+\alpha}{2},1+\alpha}([0, \infty)\times[0,h(t)])
\ \ \mbox{and }\ h\in C^1 ([0,\infty)).
\]
Moreover, we can use the classical theory for parabolic equations to obtain that $u$ is positive and
bounded in $[0,h(t))$ and $u_x(t,h(t))<0$ for $t>0$, thus $h'(t)>0$ for $t>0$. Then the following
limit is well defined
$$
\lim_{t\to\infty}h(t):=h_{\infty}\in(0,\infty].
$$
Now we list some possible situations on the asymptotic behavior of the solutions to \eqref{p}:

\begin{itemize}
	
	\item {\bf vanishing} : $h_\infty<\infty$ and $\lim_{t\to\infty}u(t,x)=0$ uniformly in $[0,h_\infty]$;
	
	\item {\bf spreading} : $h_\infty=\infty$ and $\lim_{t\to\infty}u(t,x)=1$  locally uniformly in $[0,\infty)$;
	
	\item {\bf transition} : $h_\infty=\infty$ and $\lim_{t\to\infty}|u(t,x)-U(x)|=0 \mbox{ locally uniformly in $[0,\infty)$}$,
	where $U$ is a ground state of the following elliptic equation:
\begin{equation}\label{U}
\left\{
\begin{array}{ll}
U'' + f(U)=0, & 0<x<L,\\
U'' + g(U)=0, &  x>L,\\
U'(0)=0,             \\
U( L-0)= U(L+0), \\
U'( L-0)= U'(L+0).
\end{array}
\right.
\end{equation}
\end{itemize}

\smallskip

By saying a ground state $U$ of \eqref{U}, we mean that $U$ is a positive solution to \eqref{U} satisfying
$U(\cdot)=V(\cdot-z)$ for $x>L$, with $z\in\R$ and $V$ being the unique positive symmetrically decreasing solution of
\begin{equation}\label{grs}
V''+g(V)=0\ \ \mbox{in}\ \R,\ \ V(0)=\theta^*,\ \ V'(0)=0=V(\pm \infty),
\end{equation}
where $\theta^*\in(\theta,1)$ is the constant determined by the condition
\begin{equation}\label{etta}
\int_0^{\theta^*}g(s)ds=0.
\end{equation}

We define the following three critical values:

 $$
 L_*:=\frac{1}{\sqrt{f'(0)}}\arctan\sqrt{-\frac{g'(0)}{f'(0)}}<\frac{\pi}{2\sqrt{f'(0)}}:=L_{**},
 $$
\begin{equation*}\label{L8}
L^*:=\sup\{L_0 > 0: \ \mbox{ problem \eqref{U} with $L=L_0$ has a ground state}\}.
\end{equation*}
It follows from \cite{DPS} that problem \eqref{U} has a ground state for any $0<L<L^*$ and $L^*$ is bounded.

We are now in a position to give a satisfactory explanation of the long-time dynamical behavior of problem \eqref{p}.

\begin{thm}\label{thm:dybe}
Assume that {\bf{(H)}} holds and $\phi \in \mathscr {X}(h_0)$ with $h_0>L$. Let $(u,h)$ be the solution of problem \eqref{p}
with $u_0=\sigma\phi$, and $L_*,\  L^*,\ L_{**}$ be defined as before. The following assertions hold.
\begin{itemize}
\item[\vspace{10pt}(I)] \; (Small protection zone case) Assume that $0<L\leq L_*$, then there exist $\sigma_*,\ \sigma^*\in(0,\infty)$
with $\sigma_*\leq\sigma^*$ such that the following trichotomy holds:
\begin{itemize}
\item[\vspace{10pt}(i)] \;Vanishing happens when $0<\sigma<\sigma_*;$
\item[\qquad(ii)] \;Transition happens when $\sigma\in[\sigma_*,\sigma^*];$
\item[\qquad(iii)] \;Spreading happens when $\sigma>\sigma^*$.
\end{itemize}

\item[\qquad(II)] \;(Medium-sized protection zone case) Assume that $L_*<L< \max\{L^*,\ L_{**}\}$.
\begin{itemize}
\item[\qquad(1)] \; If $L_*<L<L^*$ and $L^* < L_{**}$, then there exist $\sigma_*,\ \sigma^*\in[0,\infty)$
with $\sigma_*\leq\sigma^*$ such that the following trichotomy holds:
\begin{itemize}
\item[\vspace{10pt}(i)] \;Vanishing happens when $0<\sigma\leq \sigma_*$;
\item[\qquad(ii)] \;Transition happens when $\sigma\in(\sigma_*,\sigma^*]$;
\item[\qquad(iii)] \;Spreading happens when $\sigma>\sigma^*$.
\end{itemize}
In addition, when $h_0< R^*(L)$, then $\sigma_*>0$; when $h_0\geq R^*(L)$, then $\sigma_*=0$, where $R^*(L)$ is
given in Lemma \ref{lem:1ei2}.
\item[\qquad(2)] \; If $L^*<L<L_{**}$ and $L^* < L_{**}$, then there exists $\sigma^*\in[0,\infty)$
such that the following dichotomy holds:
\begin{itemize}
\item[\vspace{10pt}(i)] \;Vanishing happens when $0<\sigma\leq \sigma^*$;
\item[\qquad(ii)] \;Spreading happens when $\sigma>\sigma^*$.
\end{itemize}
In addition, if $h_0< R^*(L)$, then $\sigma^*>0$; while if $h_0\geq R^*(L)$, then $\sigma^*=0$.

\item[\qquad(3)] \; If $L_*<L<L_{**}$ and $L_{**}< L^*$, then there exist $\sigma_*,\ \sigma^*\in[0,\infty)$
with $\sigma_*\leq\sigma^*$ such that the following trichotomy holds:
\begin{itemize}
\item[\vspace{10pt}(i)] \;Vanishing happens when $0<\sigma\leq\sigma_*$;
\item[\qquad(ii)] \;Transition happens when $\sigma\in(\sigma_*,\sigma^*]$;
\item[\qquad(iii)] \;Spreading happens when $\sigma>\sigma^*$.
\end{itemize}
In addition, when $h_0< R^*(L)$, then $\sigma_*>0$; when $h_0\geq R^*(L)$, then $\sigma_*=0$.
\item[\qquad(4)] \; If $L_{**}<L<L^*$ and $L_{**}< L^*$, then there exists $\sigma^*\in[0,\infty)$
such that the following dichotomy holds:
\begin{itemize}
\item[\vspace{10pt}(i)] \;Transition happens when $0<\sigma\leq \sigma^*$;
\item[\qquad(ii)] \;Spreading happens when $\sigma>\sigma^*$.
\end{itemize}
\end{itemize}
\item[\qquad(III)] \;(Large protection zone case) Assume that $L>\max\{L^*,\ L_{**}\}$, then spreading happens for all $\sigma>0$.
\end{itemize}
\end{thm}

Next, we intend to study the asymptotic profiles and speeds for the solutions when spreading happens
as in Theorem \ref{thm:dybe}. The following semi-wave problem
\begin{equation}\label{eq:psemi-wave}
\left\{
 \begin{array}{ll}
 q_{c^*}''-c^*q_{c^*}'+g(q_{c^*}) =0,\ \ q_{c^*}(z)>0, & z>0,\\
 q_{c^*}(0)=0,\ \ q_{c^*}(\infty)=1, \ \ \mu q_{c^*}'(0)=c^*,
  \end{array}
 \right.
\end{equation}
will play an important role. It is well known that \eqref{eq:psemi-wave} admits
a unique solution $(c^*, q_{c^*})$ \cite{DuLou}. We call {\it $q_{c^*}$ a semi-wave with speed $c^*$}.
Based on the semi-wave, we can construct suitable upper-lower solution to deduce the following result
on the asymptotic spreading speeds and profiles of spreading solutions of \eqref{p}.
\begin{thm}\label{thm:profile of spreading sol}
Assume that {\bf{(H)}} holds and that spreading happens for a solution $(u, h)$ of problem \eqref{p}.
Let $(c^*, q_{c^*})$ be the unique solution of \eqref{eq:psemi-wave}. Then we have
\begin{equation}\label{HWt1}
\lim_{t\to\infty} h'(t) =c^*,
\end{equation}
\begin{equation}\label{WHt1}
\lim\limits_{t\to\infty} \| u(t,\cdot)- q_{c^*}(h(t) -\cdot)\| _{L^\infty ( [0, h(t)])}=0.
\end{equation}
\end{thm}

Let $\tilde{L}_*$ be given in Lemma \ref{lem:1eigenvaluemp2}. As in \cite{DPS}, it can be shown
in Section 2 that the stationary problem corresponding to \eqref{q} has a ground state provided
that $0<L<\tilde L_*$. Then
\begin{equation}\label{L8-1}
\tilde L^*:=\sup\{L_0 > 0: \ \mbox{ problem \eqref{q-statio} with $L=L_0$ has a ground state}\}
\end{equation}
is well defined and bounded.

\smallskip

Our result on problem \eqref{q} can be stated as follows.
\begin{thm}\label{thm:dybe1}
Assume that {\bf{(H)}} holds, and $u_0 \in \mathscr {X}(h_0)$ with $h_0>L_2$. Let $(u,h)$ be the solution of \eqref{q}
with $u_0=\sigma\phi$, and $\tilde{L}_*$, $\tilde{L}^*$ and $\tilde{L}_{**}$ be given in Lemma \ref{lem:1eigenvaluemp2},
\eqref{L8-1} and  Lemma \ref{lem:s1ei2}, respectively. The following assertions hold.
\begin{itemize}
\item[\vspace{10pt}(I)] \; Assume that $0<L\leq \tilde{L}_*$, then there exist $\tilde{\sigma}_*,\
\tilde{\sigma}^*\in(0,\infty)$ with $\tilde{\sigma}_*\leq\tilde{\sigma}^*$ such that the following trichotomy holds:
\begin{itemize}
\item[\vspace{10pt}(i)] \;Vanishing happens when $0<\sigma<\tilde{\sigma}_*;$
\item[\qquad(ii)] \;Transition happens when $\sigma\in[\tilde{\sigma}_*,\tilde{\sigma}^*];$
\item[\qquad(iii)] \;Spreading happens when $\sigma>\tilde{\sigma}^*$.
\end{itemize}

\item[\qquad(II)] \; Assume that $\tilde{L}_*<L< \max\{\tilde{L}^*,\ \tilde{L}_{**}\}$.
\begin{itemize}
\item[\qquad(1)] \; If $\tilde{L}_*<L<\tilde{L}^*$ and $\tilde{L}^* < \tilde{L}_{**}$, there exist
$\tilde{\sigma}_*,\ \tilde{\sigma}^*\in[0,\infty)$ with $\tilde{\sigma}_*\leq\tilde{\sigma}^*$ such that
the following trichotomy holds:
\begin{itemize}
\item[\vspace{10pt}(i)] \;Vanishing happens when $0<\sigma\leq \tilde{\sigma}_*;$
\item[\qquad(ii)] \;Transition happens when $\sigma\in(\tilde{\sigma}_*,\tilde{\sigma}^*];$
\item[\qquad(iii)] \;Spreading happens when $\sigma>\tilde{\sigma}^*$.
\end{itemize}
In addition, when $h_0< \tilde{R}^*(L)$, then $\tilde{\sigma}_*>0$; when $h_0\geq \tilde{R}^*(L)$, then
$\tilde{\sigma}_*=0$, where $\tilde{R}^*(L)$ is given in Lemma \ref{lem:s1ei2}.
\item[\qquad(2)] \; If $\tilde{L}^*<L<\tilde{L}_{**}$ and $\tilde{L}^* < \tilde{L}_{**}$, there exists
$\tilde{\sigma}^*\in[0,\infty)$ such that the following dichotomy holds:
\begin{itemize}
\item[\vspace{10pt}(i)] \;Vanishing happens when $0<\sigma\leq \tilde{\sigma}^*$;
\item[\qquad(ii)] \;Spreading happens when $\sigma>\tilde{\sigma}^*$.
\end{itemize}
In addition, if $h_0< \tilde{R}^*(L)$, then $\tilde{\sigma}^*>0$; while if $h_0\geq \tilde{R}^*(L)$,
then $\tilde{\sigma}^*=0$.

\item[\qquad(3)] \; If $\tilde{L}_*<L<\tilde{L}_{**}$ and $\tilde{L}_{**}< \tilde{L}^*$, there exist
$\tilde{\sigma}_*,\ \tilde{\sigma}^*\in[0,\infty)$ with $\tilde{\sigma}_*\leq\tilde{\sigma}^*$
such that the following trichotomy holds:
\begin{itemize}
\item[\vspace{10pt}(i)] \;Vanishing happens when $0<\sigma\leq\tilde{\sigma}_*;$
\item[\qquad(ii)] \;Transition happens when $\sigma\in(\tilde{\sigma}_*,\tilde{\sigma}^*];$
\item[\qquad(iii)] \;Spreading happens when $\sigma>\tilde{\sigma}^*$.
\end{itemize}
In addition, when $h_0< \tilde{R}^*(L)$, then $\tilde{\sigma}_*>0$; when $h_0\geq \tilde{R}^*(L)$,
then $\tilde{\sigma}_*=0$.
\item[\qquad(4)] \; If $\tilde{L}_{**}<L<\tilde{L}^*$ and $\tilde{L}_{**}< \tilde{L}^*$, there exists
$\tilde{\sigma}^*\in[0,\infty)$ such that the following dichotomy holds:
\begin{itemize}
\item[\vspace{10pt}(i)] \;Transition happens when $0<\sigma\leq \tilde{\sigma}^*$;
\item[\qquad(ii)] \;Spreading happens when $\sigma>\tilde{\sigma}^*$.
\end{itemize}
\end{itemize}
\item[\qquad(III)] \;(Large protection zone case) Assume that $L>\max\{\tilde{L}^*,\ \tilde{L}_{**}\}$,
then spreading happens for all $\sigma>0$.
\end{itemize}
Moreover, assume that spreading happens for the solution $(u, h)$ of \eqref{q}.
Let $(c^*, q_{c^*})$ be the unique solution of \eqref{eq:psemi-wave}. Then we have
\[
\lim_{t\to\infty} h'(t) =c^*\ \ \mbox{ and }\ \
\lim\limits_{t\to\infty} \| u(t,\cdot)- q_{c^*}(h(t) -\cdot)\| _{L^\infty ( [0, h(t)])}=0.
\]
\end{thm}

\smallskip

The rest of our paper is organized as follows. In Section 2, we present some preliminary results,
including the comparison principle, existence and uniqueness theorem, analysis of the associated
stationary solution problems and a general convergence result mainly due to \cite{DuLou, DM}.
In Section 3 we study the model with connected protection zone and prove Theorems \ref{thm:dybe}
and \ref{thm:profile of spreading sol}. Section 4 is devoted to studying the case of separated protection zone
and to proving Theorem \ref{thm:dybe1}. In the last section, we will give the discussion of our results
and compare our conclusions for \eqref{p} with those of other problems.

\section{Some preliminary Results}\label{sec:basic}

In this section, we present some preliminary results which will be frequently used later.

\subsection{Comparison Principle}\label{subsec:cp}
We give the following two types of comparison principle.

\begin{lem}\label{lem:comp1} Assume that $T\in (0,\infty)$, $\overline h\in C^1([0,T])$ with
$\overline h>L$ for $t\in[0,T]$, $\overline u\in C(\overline D_T)\cap C^{1,2}(D_T\backslash\{L\})$
with $D_T=\{(t,x)\in\R^2: 0<t\leq T, \ 0<x<\overline h(t)\}$, and
\begin{eqnarray*}
\left\{
\begin{array}{lll}
\overline u_{t} \geq \overline u_{xx} +f(\overline u),\; & 0<t \leq T,\
0<x<L, \\
\overline u_{t} \geq \overline u_{xx} +g(\overline u),\; & 0<t \leq T,\
L<x<\overline h(t), \\
\overline u_x(t,0)\leq 0, & 0<t \leq T,\\
\overline u(t, L-0)=\overline u(t, L+0), & 0<t \leq T,\\
\overline u_x(t, L-0)\geq\overline u_x(t, L+0), & 0<t \leq T,\\
\overline u= 0,\quad \overline h'(t)\geq -\mu \overline u_x,\quad &
0<t \leq T, \ x=\overline h(t).
\end{array} \right.
\end{eqnarray*}
If $h_0\leq \overline h(0)$, $u_0(x)\leq \overline u(0,x)$ in $[0,h_0]$,
and $(u, h)$ is a solution to \eqref{p}, then
\[
h(t)\leq\overline h(t)\mbox{ in }(0, T],\ \ u(t,x)\leq \overline u(t,x)\ \mbox{ for }t\in (0, T],\
x\in (0, h(t)).
\]
\end{lem}

\begin{lem}\label{lem:comp2} Assume that $T\in(0,\infty)$, $\overline h\in C^1([0,T])$,
$\overline u\in C(\overline D_T)\cap C^{1,2}(D_T\backslash\{L\})$ with $D_T=\{(t,x)\in\R^2:
0<t\leq T,\ r<x<\overline h(t)\}$ with $L\leq r<\overline h(t)$ for $t\in[0,T]$, and
\begin{eqnarray*}
\left\{
\begin{array}{lll}
\overline u_{t}\geq  \overline u_{xx}+ g(\overline u),\; &0<t \leq T,\ r<x<\overline h(t), \\
\overline u\geq u, &0<t \leq T, \ x= r,\\
\overline u= 0,\quad \overline h'(t)\geq -\mu \overline u_x,\quad
&0<t \leq T, \ x=\overline h(t),
\end{array} \right.
\end{eqnarray*}
with $h_0\leq \overline h(0)$, $u_0(x)\leq \overline u(0,x)$  in $[r,h_0]$.
Let  $(u, h)$ be a solution to \eqref{p}, then
\[
\mbox{ $h(t)\leq\overline h(t)$ in $(0, T]$,\quad $u(t,x)\leq
\overline u(t,x)$\ \ for $t\in (0, T]$\ and $  r<x< h(t)$.}
\]
\end{lem}

The proof of Lemma \ref{lem:comp1} is identical to that of
\cite[Lemma 5.7]{DuLin} by using  \cite[Lemmas A.2 and  A.3]{JPS};
a minor modification of such a proof yields Lemma
\ref{lem:comp2}.

\begin{remark}
\label{rem5.8}\rm The function $\overline u$, or the pair
$(\overline u, \overline h)$, in Lemmas \ref{lem:comp1}
and \ref{lem:comp2} is often called an upper solution to \eqref{p}.
A lower solution can be defined analogously by reversing all the
inequalities. The corresponding comparison principle for
lower solutions holds in each of the above cases. Similarly, we also have the corresponding
comparison principle associated with problem \eqref{q}.
\end{remark}

\subsection{Existence and uniqueness theorem}\label{sub:exist}
In this subsection, we prove that problem \eqref{p} (or \eqref{q}) admits a unique positive solution
$(u,h)$ which exists for all $t\in(0,\infty)$. Let us start with the following local existence
result.
\begin{thm}\label{thm:local}
For any given $u_0\in \mathscr {X}(h_0)$ with $h_0>L$ and any $\alpha\in (0,1)$, there is a $T>0$ such that problem
\eqref{p} $($or \eqref{q}$)$ admits a solution $(u,h)$ which satisfies
$$u\in C^{\frac{1+\alpha}{2},1+\alpha}([0,T]\times[0,h(t)])\cap C^{1,2}((0,T]\times([0,h(t)]\backslash\{L\})),\ \ h\in C^{1+\frac{\alpha}{2}}([0,T]).$$
\end{thm}

\begin{proof} We only treat system \eqref{p}; the analysis for system \eqref{q} is similar.
In fact, our argument mainly follows that of \cite[Theorem 2.1]{DuLin} but with necessary modifications.
We divide the proof into two steps as follows.

\smallskip

$Step\ 1$. As in \cite{DuLin}, we first  straighten the free boundary. Denote $\delta:=h_0-L$ and let
$\xi(y)$ be a nonnegative function in $C^{3}(\R)$ such that
\[
\xi(y)=1\  \mbox{ if }\  | y-h_0|\leq \frac{\delta}{4},\  \xi(y)=0\ \mbox{ if } \ |y-h_0| \geq \frac{\delta}{2},\ \
\mbox{and}\ \ |\xi'(y)|<\frac{6}{\delta}\ \mbox{ for }\  y\in \R.
\]
For any given $h(t)$, we can straighten the free boundary in \eqref{p} by the transformation $(t,y)\to (t,x)$, where
\[
x=y+\xi(y)(h(t)-h_0), \ \ y\geq 0.
\]
We can choose the small  $t\in(0, \mathcal{T}]$ such that
\[
|h(t)-h_0|\leq \frac{\delta}{8}.
\]
This implies that the above transformation is a diffeomorphism from $[0,\infty)$ onto $[0,\infty)$, and
\begin{align*}
 0\leq x\leq L+\frac{\delta}{2} &\Leftrightarrow  0\leq y\leq L+\frac{\delta}{2},\\
 L+\frac{\delta}{2}\leq x \leq h(t) &\Leftrightarrow  L+\frac{\delta}{2}\leq y\leq h_0,\\
 x=h(t) &\Leftrightarrow y=h_0.
\end{align*}
Furthermore, some direct calculations give
\[
\frac{\partial y}{\partial x}=\frac{1}{1+\xi'(y)(h(t)-h_0)}\equiv \sqrt{A(h(t),y)},
\]
\[
\frac{\partial^2 y}{\partial x^2}=-\frac{\xi''(y)(h(t)-h_0)}{[1+\xi'(y)(h(t)-h_0)]^3}\equiv B(h(t),y),
\]
\[
-\frac{1}{h'(t)}\frac{\partial y}{\partial t}=\frac{\xi(y)}{1+\xi'(y)(h(t)-h_0)}\equiv C(h(t),y).
\]
It follows that
\begin{equation}\label{eq:abx1}
A(h(t),y)\equiv 1\ \ \mbox{ and }\ \ B(h(t),y)\equiv0\equiv C(h(t),y) \ \
\mbox{for }\ t>0, 0\leq y\leq L+\frac{\delta}{2}.
\end{equation}
If we set
\[
w(t,y):=u(t,y+\xi(y)(h(t)-h_0))=u(t,x),
\]
then the free boundary problem \eqref{p} becomes
\begin{equation}\label{lin1}
\left\{
\begin{array}{ll}
 w_t -\tilde{A}(t,y)w_{yy} - \tilde{B}(t,y)w_{y}=f(w(t,y)), & t\in(0,\mathcal{T}], y\in(0,L), \\
 w_t -\tilde{A}(t,y)w_{yy} - \tilde{B}(t,y)w_{y}=g(w(t,y)), &t\in(0,\mathcal{T}],  y\in(L,h_0), \\
 w_y(t,0)=w(t,h_0)=0,  &  t\in(0,\mathcal{T}],\\
 h'(t) = -\mu w_y(t,h_0), & t\in(0,\mathcal{T}],\\
 w(t,L-0)=w(t,L+0), & t\in(0,\mathcal{T}],\\
 w_y(t,L-0)=w_y(t,L+0), & t\in(0,\mathcal{T}],\\
 h(0)=h_0,\ \  w(0,y) =u_0(y), & y\in[0,h_0],
\end{array}
\right.
\end{equation}
where $\tilde{A}(t,y):=A(h(t),y)$, $\tilde{B}(t,y):=B(h(t),y)+h'(t)C(h(t),y)$ and \eqref{eq:abx1} are used.

Denote $h_1:=-\mu u'_{0}(h_0)$ and $\mathcal{T}:=\frac{\delta}{16(1+ h_1)}$.
For $T\in(0, \mathcal{T}]$, we define $\Omega_{T}:=[0,T]\times[0,h_{0}]$,
\begin{align*}
& D^{w}_{T}:=\{w\in C(\Omega_{T}):\ w(0,y)=u_{0}(y),\ \| w-u_{0}\|_{C(\Omega_{T})}\leq1\},\\
& D^{h}_{T}:=\{h\in C^{1}([0,T]):\ h(0)=h_{0},\ h'(0)=h_{1},\ \| h'-h_{1}\|_{C([0,T])}\leq1\}.
\end{align*}
Then we can notice that $D_{T}=D^{w}_{T}\times D^{h}_{T}$ is a complete metric space with the  metric
\[
d((w,h),(\tilde{w},\tilde{h}))=\|w-\tilde{w}\|_{C(\Omega_{T})}+\|h'-\tilde{h}'\|_{C([0,T])}.
\]

For any given $(w,h)\in D_{T}$ with $T\in(0,\mathcal{T}]$, we extend $(w,h)$ to $t>T$ by defining
\[
(w(t,y),h(t))=(w(T,y),h(T))\ \mbox{ for }\ t>T, \ y\in[0,h_0],
\]
and we extend the associated $\tilde{A}(t,y)$ and $\tilde{B}(t,y)$ similarly. For simplicity the extended functions are still
denoted by themselves.

Let $L_f$ and $L_g$ be the Lipschitz constants of $f$ and $g$, respectively. For the fixed $T\in(0,\mathcal{T}]$, we choose $(w,h)\in D_{T}$ to consider the following problem:
\begin{equation}\label{1n3}
\left\{
\begin{array}{ll}
\tilde{w}_t -\tilde{A}(t,y)\tilde{w}_{yy} - \tilde{B}(t,y)\tilde{w}_{y}=f(w(t,y)), &  t\in(0,\mathcal{T}],\ y\in(0,L),\\
 \tilde{w}_t -\tilde{A}(t,y)\tilde{w}_{yy} - \tilde{B}(t,y)\tilde{w}_{y}=g(w(t,y)), &  t\in(0,\mathcal{T}],\ y\in(L,h_0),\\
 \tilde{w}_y(t,0)=\tilde{w}(t,h_0)=0,  &  t\in(0,\mathcal{T}],\\
 \tilde{w}(t,L-0)=\tilde{w}(t,L+0), & t\in(0,\mathcal{T}],\\
 \tilde{w}_y(t,L-0)=\tilde{w}_y(t,L+0), & t\in(0,\mathcal{T}],\\
 \tilde{w}(0,y) =u_0(y), & y\in[0,h_0],
\end{array}
\right.
\end{equation}
where the above extensions of $(w,h)$, $\tilde{A}$ and $\tilde{B}$ are assumed for $t>T$. We apply the standard $L^p$ theory to find that \eqref{1n3} admits a unique solution $\tilde{w}(t,y)\in W_p^{1,2}([0,\mathcal{T}]\times(0,h_0))$ for any $p>1$, and there exists
$C_0>0$ depending only on $p$, $\mathcal{T}$, $h_0$, $\| u_{0}\|_{C^2([0,h_{0}])}$, $L_f$ and $L_g$
such that
\begin{equation}\label{20200915}
\|\tilde{w}(t,y)\|_{W_p^{1,2}([0,\mathcal{T}]\times(0,h_0))}\leq C_0.
\end{equation}
In view of the Sobolev embedding theorem, for any $\alpha\in(0,1)$, there exists $K>0$ depending on $p>1$, $h_0$, $\alpha$ and $\mathcal{T}$ such that
\[
 \|\tilde{w}\| _{C^{\frac{1+\alpha}{2},{1+\alpha}}(\Omega_{\mathcal{T}})}\leq K \|\tilde{w}(t,y)\|_{W_p^{1,2}([0,\mathcal{T}]\times(0,h_0)))}\leq KC_0=:C_{1}.
\]
It follows that
\begin{equation}\label{eq1}
 \|\tilde{w}\| _{C^{\frac{1+\alpha}{2},{1+\alpha}}(\Omega_{T})}\leq C_{1}.
\end{equation}
For any $t\in[0,\mathcal{T}]$, define
\[
\tilde{h}(t)=h_0-\int_0^t \mu\tilde{w}_{y}(\tau, h_{0})d\tau,
\]
we deduce
\[
\tilde{h}'(t)=-\mu\tilde{w}_{y}(t, h_{0}),\ \ \tilde{h}(0)=h_0,\ \ \tilde{h}'(0)=-\mu\tilde{w}_{y}(0, h_{0})=h_1,
\]
and thus $\tilde{h}'\in C^{\frac{\alpha}{2}}([0,\mathcal{T}])$, which satisfies
\begin{equation}\label{eq2}
\|\tilde{h}'\| _{C^{\frac{\alpha}{2}}([0,\mathcal{T}])}\leq\mu C_{1}=:C_{2}.
\end{equation}

\smallskip

$Step\ 2$.  We define $\tilde{\mathcal{F}}: D_T\to C(\Omega_{\mathcal{T}})\times C^1([0,\mathcal{T}])$ and
$\mathcal{F}: D_T\to C(\Omega_{T})\times C^1([0,T])$ by
\[
\tilde{\mathcal{F}}(w,h)=(\tilde{w},\tilde{h}),\ \ \ \mathcal{F}(w,h)=(\tilde{w},\tilde{h})|_{\{t\in[0,T]\}}.
\]
Clearly $(w,h)$ is a solution of \eqref{lin1} for $t\in[0,T]$ if and
only if it is a fixed point of $\mathcal{F}$ in $D_{T}$. We will show that $\mathcal{F}$ maps $D_{T}$ into itself provided that
$ T>0$ is small enough.

By \eqref{eq1} and \eqref{eq2}, we have
\[
\|\tilde{h}'-h_{1}\|_{C([0,T])}\leq T^{\frac{\alpha}{2}}\|\tilde{h}'\| _{C^{\frac{\alpha}{2}}([0,T])}\leq T^{\frac{\alpha}{2}}C_{2}
\]
and
\[
\|\tilde{w}-u_{0}\|_{C(\Omega_{T})}\leq T^{\frac{1+\alpha}{2}} (C_{1}+\tilde{C}_1\|u_0\|_{C^2([0,h_0])})=:T^{\frac{1+\alpha}{2}}\hat{C},
\]
where $\tilde{C}_1$ depends on $h_0$ and $\alpha$. We take $T\leq T_{1}=\min\{\ \hat{C}^{-\frac{2}{\alpha}}_{2},\ C^{-\frac{2}{1+\alpha}}_{1},\ \mathcal{T}\}$,
then we can show that $\mathcal{F}$ maps $\Omega_{T}$ into itself.

Next, we prove that $\mathcal{F}$ is a contraction mapping on $\Omega_{T}$
for $T>0$ small enough. In fact, choose $(w_i, h_i)\in D_{T}$ for $i=1,\ 2$
and denote
\[
(\tilde{w}_{i}, \tilde{h}_{i})=\tilde{\mathcal{F}}(w_{i}, h_{i}).
\]
We assume that $(w_i, h_i)$ are extended to $t>T$ as before. We denote the associated $\tilde{A}(t,y)$ and $\tilde{B}(t,y)$ by
$\tilde{A}_i(t,y)$ and $\tilde{B}_i(t,y)$  and assume that they are also extended to $t>T$ as before.

Denote $W(t,y):=\tilde{w}_{1}(t,y)-\tilde{w}_{2}(t,y)$, then we find that $W(t,y)$ satisfies that
\begin{equation}\label{p1}
\left\{
\begin{array}{ll}
 W_{t} -\tilde{A}_2(t,y)W_{yy} - \tilde{B}_2(t,y)W_{y}=\tilde{f}(t,y), & t\in(0,\mathcal{T}],\, y\in(0, L), \\
 W_{t} -\tilde{A}_2(t,y)W_{yy} - \tilde{B}_2(t,y)W_{y}=\tilde{g}(t,y), & t\in(0,\mathcal{T}],\ y\in(L, h_0), \\
 W_y(t,0)=W(t,h_0)= 0,  &  t\in(0,\mathcal{T}],\\
 W(t,L-0) =W(t,L+0), & t\in(0,\mathcal{T}],\\
 W_y(t,L-0) =W_y(t,L+0), & t\in(0,\mathcal{T}],\\
 W(0,y)=0,  & y\in[0, h_0],
\end{array}
\right.
\end{equation}
where
\[
\tilde{f}(t,y)=[\tilde{A}_1(t,y)-\tilde{A}_2(t,y)](\tilde{w}_{1})_{yy}+ [\tilde{B}_1(t,y)-\tilde{B}_2(t,y)](\tilde{w}_1)_y+f(w_1(t,y))-f(w_2(t,y)),
\]
and
\[
\tilde{g}(t,y)=[\tilde{A}_1(t,y)-\tilde{A}_2(t,y)](\tilde{w}_{1})_{yy}+ [\tilde{B}_1(t,y)-\tilde{B}_2(t,y)](\tilde{w}_1)_y+g(w_1(t,y))-g(w_2(t,y)).
\]
Due to \eqref{20200915}, we use the $L^{p}$  ( for any $p>1$) estimate to obtain that
\[
\| W\|_{w_{p}^{1,2}((0,\mathcal{T}]\times(0,h_0))}\leq C_3 (\| w_{1}-w_{2}\|_{C(\Omega_{T})}+\| h'_{1}-h'_{2}\|_{C([0,T])})
\]
with $C_3$ depending only on $\mathcal{T},\ h_{0}$, $\|u\|_{C^2([0,h_0])}$, $L_f$, $L_g$ and $p$.
By Sobolev embedding theorem, we obtain that
\[
\| W\|_{C^{\frac{1+\alpha}{2},{1+\alpha}}(\Omega_{\mathcal{T}})}\leq C_4 (\| w_{1}-w_{2}\|_{C(\Omega_{T})}+\| h'_{1}-h'_{2}\|_{C([0,T])})
\]
with $C_4$ depending only on $\mathcal{T},\ h_{0}$, $C_3$ and $\alpha$.
Thus for any $T\in(0, T_1]$, we deduce that
\[
\| W\|_{C^{\frac{1+\alpha}{2},{1+\alpha}}(\Omega_{T})}\leq C_4 (\| w_{1}-w_{2}\|_{C(\Omega_{T})}+\| h'_{1}-h'_{2}\|_{C([0,T])}).
\]
It follows that
\begin{align*}
\| \tilde{w}_{1}- \tilde{w}_{2}\|_{C(\Omega_{T})}&\leq T^{\frac{1+\alpha}{2}}\| \tilde{w}_{1}-\tilde{w}_{2}\|_{C^{\frac{1+\alpha}{2},{1+\alpha}}(\Omega_{T})}\\
 &\leq C_4 T^{\frac{1+\alpha}{2}} (\| w_{1}-w_{2}\|_{C(\Omega_{T})}+\| h'_{1}-h'_{2}\|_{C([0,T])}),
\end{align*}
and
\begin{align*}
\| \tilde{h}'_{1}-\tilde{h}'_{2}\|_{C ([0,T])}&\leq \mu (\| (\tilde{w}_{1})_{y}- (\tilde{w}_{2})_{y}\|_{C(\Omega_{T})})\\
&\leq\mu T^{\frac{\alpha}{2}} \|\tilde{w}_{1}-\tilde{w}_{2}\|_{C^{\frac{1+\alpha}{2},{1+\alpha}}(\Omega_{T})}\\
 &\leq C_4\mu T^{\frac{\alpha}{2}} (\| w_{1}-w_{2}\|_{C(\Omega_{T})}+\| h'_{1}-h'_{2}\|_{C([0,T])}).
\end{align*}
Therefore, we can obtain that
\begin{align*}
d((\tilde{w}_{1},\tilde{h}_{1}),(\tilde{w}_{2},\tilde{h}_{2}))&=\| \tilde{w}_{1}- \tilde{w}_{2}\|_{C(\Omega_{T})}+\| \tilde{h}'_{1}-\tilde{h}'_{2}\|_{C ([0,T])}\\
& \leq (\mu T^{\frac{\alpha}{2}}+ T^{\frac{1+\alpha}{2}})C_4 (\| w_{1}-w_{2}\|_{C(\Omega_{T})}+\| h'_{1}-h'_{2}\|_{C([0,T])})\\
&= (\mu T^{\frac{\alpha}{2}}+ T^{\frac{1+\alpha}{2}})C_4 d((w_{1},h_{1}),(w_{2},h_{2})).
\end{align*}
Set
\[
T_{2}:=\min\Big\{1,\ T_{1},\ \Big(\frac{1}{2(\mu+1) C_4}\Big)^{\frac{2}{\alpha}},\ \Big(\frac{1}{2(\mu+1)C_4}\Big)^{\frac{2}{1+\alpha}},\ \frac{\delta}{16(1+ h_1)}\Big\},
 \]
then for any $T\in(0,T_{2}]$, we have
 \[
 d((\tilde{w}_{1},\tilde{h}_{1}),(\tilde{w}_{2},\tilde{h}_{2}))\leq \frac{1}{2} d((w_{1},h_{1}),(w_{2},h_{2})),
 \]
which implies that $\mathcal{F}$ is a contraction mapping on $D_T$ for any $T\in(0,T_{2}]$. It follows from the contraction mapping theorem that $\mathcal{F}$
has a unique fixed point $(w, h)\in D_{T}$, which gives a unique solution to \eqref{lin1}. Moreover,  by the Schauder estimates, we have additional
regularity for $(w,h)$ as a solution of \eqref{lin1}, namely, $h\in C^{1+\alpha/2}((0,T])$ and $w\in C^{1,2}([0,T]\times([0,h_0]/\{L\}))$, and
\eqref{eq1} and \eqref{eq2} hold. In other words, we can get a unique local classical solution of \eqref{p} through the transformation.
\end{proof}

\begin{lem}\label{lem:global}
Problem \eqref{p} $($or \eqref{q}$)$ admits a unique positive solution $(u, h)$ which exists for
all $t\in(0, \infty)$.
\end{lem}

\begin{proof}
We only prove the assertion for \eqref{p}, and the analysis for system \eqref{q} is similar. Let $[0, T_{max})$
be the maximal time interval of the existence of the solution. In view of Theorem \ref{thm:local}, it remains
to show $T_{max}=\infty$.

Suppose for contradiction that $T_{max}<\infty$.
Set $K_0:=\|u_0\|_\infty+1$. Since $f(s)<0$ and $g(s) <0$ for $s>1$, then the comparison principle (Lemma \ref{lem:comp1}) gives
\[
0\leq u(t,x)\leq K_0 \ \mbox{ for } t\in(0,T_{max}), \ x\in (0,h(t)).
\]
Since $g$ is globally Lipschitz and $g(0)=0$, there exists $K_1>0$ depending on $K_0$
such that $g(s)\leq K_1$ for $s\in[0, K_0]$. Let us construct the auxiliary function
\[
v(t,x)=K_0\big[2M(h(t)-x)-M^{2}(h(t)-x)^{2}\big]
\]
over the region
\[
\Omega_{M}=\{ (t,x):\ 0<t<T_{max},\ h(t)-M^{-1}<x<h(t)\}
\]
with $M>0$ to be determined.

Firstly, we choose $M\geq \frac{1}{h_0-L}$. Then  $h(t)-M^{-1}\geq L$ for $t\geq 0$.
Since $0\leq v(t,x)\leq K_0$ in $\Omega_{M}$, a direct calculation yields
\[
v_{t}-v_{xx}\geq2K_0M^{2}\geq K_1\geq g(v)\ \mbox{ in } \Omega_{M},
\]
if $M\geq\sqrt{\frac{K_1}{2K_0}}$. On the other hand, $v(t,h(t))=0=u(t,h(t))$, and
\[
v(t,h(t)-M^{-1})=K_0\geq u(t,h(t)-M^{-1}).
\]
Since $v(0,h_0)=0=u_0(h_0)$, if $M\geq \frac{4\| u_{0}\|_{C^{1}([0,h_0])}}{3K_0}$,
then one can check that
\[
v(0,x)\geq u_0(x)\ \ \mbox{ for } x\in[h_{0}-M^{-1},h_{0}].
\]
Therefore by choosing
\[
M:=\max\Big\{\frac{1}{h_0-L},\ \sqrt{\frac{K_1}{2K_0}},\ \frac{4\| u_{0}\|_{C^{1}}}{3K_0}\Big\},
\]
and using Lemma \ref{lem:comp2}, we can deduce that $v(t,x)\geq u(t,x)$ in $\Omega_{M}$. It then follows that
\[
h'(t)=-\mu u_{x}(t,h(t))\leq-\mu v_{x}(t,h(t))= 2\mu M K_0:= K_2.
\]

Let us  now fix $\epsilon\in(0,T_{max})$. Similar to the proof of Theorem \ref{thm:local}, by the $L^p$ and Schauder estimates together
with the Sobolev embedding theorem for parabolic equation, we can find $C_1>0$ depending only on $\epsilon$, $T_{max}$, $h_0$,
$\| u_0\|_{C^1}$ and $K_2$ such that
\[
||u||_{C^{1,2}([\varepsilon,T_{max})\times([0, h(t)]/\{L\}))}\leq C_1.
\]
This implies that $(u, h)$ exists on $[0,T_{max})$. Choosing $\{t_n\}^{\infty}_{n=1}\subseteq(0,T_{max})$ with $t_n\nearrow T_{max}$, and regarding
$t_n$ and $(u(t_n, x), h(t_n))$  $(n\geq 1)$ as the initial time and initial datum, respectively, it then follows from
the proof of Theorem \ref{thm:local} that there exists $s_0>0$ depending on $K_0$, $K_2$ and $C_1$ independent of $n$ such that
problem \eqref{p} has a unique solution $(u, h)$ in $[t_n, t_n+s_0)$ for all $n\geq 1$. This indicates that the solution $(u, h)$ of \eqref{p}
can be extended uniquely to $[0,t_n+s_0)$. Hence $t_n+s_0>T_{max}$ when $n$ is large, which contradicts the definition of $T_{max}$.
The proof of this lemma is thus complete.
\end{proof}

\subsection{Stationary solutions}\label{subsect:statin}
A stationary solution of \eqref{p} is a solution of \eqref{U}. By a similar argument as in
\cite[Lemma 2.2]{DPS}, we obtain the following result.
\begin{lem}\label{lem:staso}
Assume that {\bf{(H)}} holds. For any $L>0$, all solutions of the stationary problem \eqref{U} are one of the following types:
\begin{itemize}
\item[\vspace{10pt}(1)] \;Trivial solution: $U\equiv 0$;
\item[\qquad(2)] \;Positive constant solution: $U\equiv 1$;
\item[\qquad(3)] \;Ground states: $U$ is positive and decreasing for $x\in[0,\infty)$, and
\[\left\{
\begin{array}{ll}
 U'' + f(U)=0, & 0<x<L,\\
 U'' + g(U)=0, &  x>L,\\
 U'(0)= 0, \\
 U(L-0)= U(L+0), \\
 U'( L-0)= U'(L+0).
\end{array}
\right.
\]
Moreover, when $x>L$, $U(\cdot)=V(\cdot-z)$ where $z\in\R$ and $V$ is the unique
positive symmetrically decreasing solution of \eqref{grs}.
\item[\qquad(4)] \;Periodic solutions: $U(x)>0$ for $x\geq 0$ and
$U(x)=P(x-\hat{z})$ for $x>L$, where $\hat{z}\in\R$, $P$ is a periodic solution of
$P'' + g(P)=0$ with  $0<\min P<\theta <\max P<\theta^*$.
\end{itemize}
\end{lem}

When the stationary solution $U(x)$ is a ground state, it then follows that
\begin{equation*}\label{Utot}
\|U\|_{L^\infty([0,\infty))}=U(0)< \theta^*.
\end{equation*}

A ground state $U$ of \eqref{q} is a positive solution of
the following elliptic problem:
\begin{equation}\label{q-statio}
\left\{
\begin{array}{ll}
 U'' + f(U)=0, &   x \in (L_1, L_2),\\
 U'' + g(U)=0, &  x\in(0, L_1)\cup(L_2,\infty),\\
 U'(0)=0, & \\
 U( L_i-0)= U( L_i+0), & i=1, 2, \\
 U'( L_i-0)= U'( L_i+0), & i=1, 2,
\end{array}
\right.
\end{equation}
when $x>L_2$, $U(\cdot)=V(\cdot-z_1)$, where $z_1\in\R$ and $V$ is the unique positive decreasing solution
of \eqref{grs}, and when $x\in[0,L_1]$, $U(\cdot)=P(\cdot-z_2)$ where $z_2\in\R$ and $P$ is a periodic
solution of $P'' + g(P)=0$ satisfying  $0<\min P<\theta <\max P<\theta^*$.

Indeed, \eqref{q-statio} may have eight types of ground state; one may refer to \cite{DPS} for more details.
Any ground state $U$ of \eqref{q-statio} satisfies
\[\|U\|_{L^\infty([0,\infty))}< \theta^*.\]

\subsection{A general convergence theorem}\label{subsect:zeronumber}

Following a similar analysis as in \cite{ DuLou,DM}, we are able to state the following general convergence result.

\begin{thm}\label{thm:convergence}{\rm(Convergence theorem for systems \eqref{p} and \eqref{q})}
Assume that $(u,h)$ is a solution of \eqref{p} $($or \eqref{q}$)$ with $u_0 \in \mathscr{X}(h_0)$ and
$h_0>L$ $($or $h_0>L_2$$)$. Then $u$ converges to a solution $U$ of \eqref{U} $($or \eqref{q-statio}$)$ as $t\to \infty$
locally uniformly in $[0,h_\infty)$ with $h_\infty<\infty$ or $h_\infty=\infty$, where $U$ is one
of the following types: $0$, $1$ and  ground states of \eqref{U} $($or \eqref{q-statio}$)$. Moreover,
if $h_\infty<\infty$, then $U\equiv 0$.
\end{thm}

\begin{proof}\ We only sketch the proof for system \eqref{p}; the analysis for system \eqref{q} is similar.

Let $\omega(u)$ be the $\omega$-limit set of $u(t,\cdot)$ in the topology of $L^{\infty}_{loc}([0,h_\infty))$.
By the local parabolic estimates, the definition of $\omega(u)$ remains unchanged if the topology of
$L^{\infty}_{loc}([0,h_\infty))$ is replaced by that of $C^2_{loc}([0,L)\cup(L,h_\infty))\cap C^1_{loc}([0,h_\infty))$.
Clearly, $\omega(u)$ is a compact, connected and invariant set.
By the argument of \cite[Theorem 1.1]{DM} and \cite[Theorem 1.1]{DuLou} with slight modifications, it can be shown that
$\omega(u)$ consists of only one element, which is either a constant solution or a decreasing solution of \eqref{U}.
In view of Lemma \ref{lem:staso}, $\omega(u)$ contains either $0$, or $1$, or a ground state of \eqref{U}.

Finally, we claim that if $h_\infty<\infty$, $\omega(u)=\{0\}$. Otherwise, $\omega(u)$ contains a nontrivial nonnegative solution $v$ of the problem
\begin{equation}\nonumber
 \left\{
 \begin{array}{ll}
 v_{xx}+f(v)=0, & 0< x <L,\\
 v_{xx}+g(v)=0, & L<x<h_\infty,\\
 v_x(0)=v(h_\infty)=0,\\
 v(L-0)=v(L+0),\\
 v_x(L-0)=v_x(L+0).
 \end{array}
 \right.
 \end{equation}
Due to $f(0)=g(0)=0$, it follows from the strong maximum principle and the Hopf boundary lemma that
$v>0$ in $[0,h_\infty)$ and $v'(h_\infty)<0$. By the definition of $\omega(u)$, we know that
$u(t,\cdot)\to v(\cdot)$ in $C^1_{loc}([0,h_\infty))$ as $t\to \infty$. Using a similar argument as in \cite[Theorem 2.1]{DuLin}
by straightening the free boundary one can show that
$$
\|u(t,\cdot)-v(\cdot)\|_{C^1((0, h(t)])}\to 0\ \ \mbox{as } t\to\infty.
$$
It follows that
\[
\lim_{t\to\infty}h'(t)=-\mu \lim_{t\to\infty}u_x(t,h(t))=-\mu v'(h_\infty)=\delta
\]
for some $\delta>0$. This contradicts the assumption that $h_\infty<\infty$, which ends the proof.
\end{proof}

\section{A Connected Protection Zone: Proof of Theorems \ref{thm:dybe} and \ref{thm:profile of spreading sol}} \label{sec:spanv}

In this section, we treat the case of connected protection zone and prove Theorems \ref{thm:dybe} and \ref{thm:profile of spreading sol}.

\subsection{Sufficient conditions for spreading}\label{subs:spr}
By the phase plane analysis as in \cite[Subsection 2.1]{DPS}, we obtain that for any $\alpha\in(\theta^*,1)$ with $\theta^*$
be given in \eqref{etta}, the following problem
\begin{equation}\label{vala}
v_\alpha''+g(v_\alpha)=0<v_\alpha\leq \alpha\ \ \mbox{ in } (0,2l_\alpha),
\ \ \  v_\alpha(0)=v_\alpha(2l_\alpha)=0,\ \ \ v_\alpha(l_\alpha)=\alpha,
\end{equation}
with
\begin{equation}\label{la1}
l_\alpha=\int_0^\alpha\frac{ds}{\sqrt{2\int_s^\alpha g(v)dv}}\in(0,\infty),
\end{equation}
admits a solution $v_\alpha$, which is used to construct suitable lower solutions
to give the following sufficient condition for spreading.
\begin{lem}\label{lem:spreadc}
If for each $\alpha\in(\theta^*,1]$, $u_0(x)\geq v_\alpha(x-r)$ in $[r,r+2l_\alpha]$ and
$h_0\geq r+2l_\alpha$ for some $r\geq L$, where $v_\alpha$ is a solution of \eqref{vala} and
$l_\alpha$ is given in \eqref{la1}, then spreading happens for $(u, h)$.
\end{lem}
\begin{proof}
It follows from the comparison principle that $u(t,x)\geq v_\alpha(x-r)$ for $x\in[r,r+2l_\alpha]$ and
$t>0$. As $v_\alpha(l_\alpha)=\alpha\in(\theta^*,1]$, only $1$ is possible solution of \eqref{U} bigger than
$v_\alpha(x-r)$ for $x\in[r,r+2l_\alpha]$, the conclusion follows from Theorem \ref{thm:convergence} immediately.
\end{proof}

Our second sufficient condition ensuring spreading is that the initial datum is sufficiently large
on any given interval.

\begin{lem}\label{lem:spread2}
Assume that $h_0>L$ and $\phi\in \mathscr{X}(h_0)$. Let $(u, h)$ be the solution of \eqref{p} with
initial datum $u_0=\sigma\phi$. Then spreading happens provided that $\sigma$ is sufficiently large.
\end{lem}
\begin{proof}
 It is known that the following eigenvalue problem
\begin{equation*}\label{SL-half}
\left\{
 \begin{array}{l}
 -\psi''(x) -\frac{1}{2}\psi'(x)=\kappa \psi (x) ,\quad x\in (\frac{1}{2},1),\\
 \psi' (\frac{1}{2}) =\psi (1)=0,
 \end{array}
\right.
\end{equation*}
has the principal eigenvalue $\kappa_1$ and the corresponding principal eigenfunction $\psi_1$ satisfies
$\psi_1(x)>0$ and $\psi'_1 (x) <0$ for all $x\in (\frac{1}{2},1]$. We assume
further that
\[
\|\psi_1\|_{L^\infty([\frac{1}{2},1])} =\psi_1 \big(\frac{1}{2}\big)=1.
\]
Let us extend $\psi_1$ to $[0,1]$ with $\psi_1(x)=\psi_1(1-x)$ for $x\in[0,\frac{1}{2}]$. Then clearly
\begin{equation*}\label{SL}
\left\{
 \begin{array}{l}
 \psi_1''(x) +\frac{{\rm sgn}\big(x-\frac{1}{2}\big)}{2}\psi_1'(x) +\kappa_1
 \psi_1 (x) =0,\quad x\in (0,1),\\
 \psi_1' (\frac{1}{2})=\psi_1 (0) =\psi_1 (1)=0.
 \end{array}
\right.
\end{equation*}

For any given $\alpha\in(\theta^*,1)$, let $v_\alpha$ be a solution of \eqref{vala} and $l_\alpha$ be given in
\eqref{la1}. As $g$ is globally Lipschitz on $[0,\infty)$, we can find $M>0$ such that $g(u)\geq -Mu$ for all $u\geq 0$.

Choose positive constants $\varepsilon$, $T$, $\gamma$ and $\rho$ as follows:
\[
 \varepsilon:= \frac{1}{2}\min \big\{1,\ (h_0-L)^2\big\},\; \ T:= (L+2l_\alpha)^2,\;\ \gamma=|\kappa_1|+M(T+1),
\]
  and
\begin{equation}\label{def:rho-1}
\frac{\rho}{(T+\varepsilon)^\gamma}\psi_1\Big(\frac{x}{\sqrt{T+\varepsilon}}\Big)\geq
v_\alpha(x)\ \ \mbox{in } [0,2l_\alpha],\ \ - 2\mu \rho \psi'_1 (1) > (T+1)^\gamma.
\end{equation}
Define
\begin{equation*}\label{lower-sol}
w(t,x) := \frac{\rho}{(t+\varepsilon)^\gamma} \ \psi_1 \left(
\frac{x-L}{\sqrt{t+\varepsilon}}\right) \quad \mbox{for } t\geq 0,\ x\in
[L,L+\sqrt{t+\varepsilon}].
\end{equation*}
For $t\in [0,T]$ and $x\in (L,L+\sqrt{t+\varepsilon})$,
 we have
\begin{eqnarray*}
w_t -w_{xx} - g(w) & \leq & w_t -w_{xx} +M w \\
& = & \frac{-\rho}{(t+\varepsilon)^{\gamma+1}} \left[ \psi''_1
+\frac{x-L}{2\sqrt{t+\varepsilon}} \psi'_1  +(\gamma -
M(t+\varepsilon)) \psi_1 \right]\\
& \leq & \frac{-\rho}{(t+\varepsilon)^{\gamma+1}} \left[
\psi''_1 +\frac{{\rm sgn}(\frac{x-L}{\sqrt{t+\varepsilon}}-\frac{1}{2})}{2} \psi'_1  + \kappa_1
\psi_1 \right] =0.
\end{eqnarray*}
Clearly $w(t,L)=w(t,L+\sqrt{t+\varepsilon})=0$ for $t>0$, and by \eqref{def:rho-1}
we have that for $t\in[0,T]$,
$$
(L+\sqrt{t+\varepsilon})' + \mu w_x (t, L+\sqrt{t+\varepsilon})
\leq \frac{1}{2\sqrt{t+\varepsilon}} \left[ 1 + \frac{2\mu
\rho}{(T+1)^\gamma } \psi'_1 (1) \right] <0.
$$
Since $\varepsilon <(h_0-L)^2$ we can choose $\sigma_1>0$ large such that
$$
w(0,x) = \frac{\rho}{\varepsilon^\gamma} \psi_1 \Big(
\frac{x-L}{\sqrt{\varepsilon}}\Big) < \sigma_1\phi(x)\quad
\mbox{for } x\in [L, L+\sqrt{\varepsilon}]\subset
[0, h_0].
$$

Hence $(w(t,x), L+\sqrt{t+\varepsilon})$ is a lower solution of \eqref{p} over $[0,T]\times[L,L+\sqrt{t+\varepsilon}]$
and the comparison principle implies that
\[
u(T,x;\sigma_1\phi)\geq w(T,x)\geq v_\alpha(x-L) \ \mbox{ for } x\in[L,L+2l_\alpha]\subset[L, h(T;\sigma_1\phi)).
\]
This and Lemma \ref{lem:spreadc} yield that spreading happens for this $\sigma_1$. The lemma is proved.
\end{proof}
Assume that $h_0 >L$ and $\phi \in \mathscr{X}(h_0)$. Let $u_0=\sigma\phi$ with $\sigma>0$. Then
\[
\Sigma_1=\big\{\sigma_1>0: \mbox{ spreading happens for } \sigma\in[\sigma_1,\infty)\}
\]
is well defined.
We end this subsection with the following property of the set $\Sigma_1$.
\begin{lem}\label{lem:sd112}
For any given $L>0$, $\Sigma_1$ is a nonempty open interval.
\end{lem}
\begin{proof}
It follows from Lemma \ref{lem:spread2} that spreading happens for  large $\sigma>0$, then $\Sigma_1$ is non-empty.

We next prove that $\Sigma_1$ is an open interval. If $\sigma_1\in \Sigma_1$, there is $T_1>0$ large enough
such that for any given $\alpha\in(\theta^*,1)$,
\begin{equation}
\label{u-v_Z}
 u_1(T_1, x)>v_\alpha(x) \mbox{ in } [L,L+2l_\alpha]\ \ \mbox{ and }\ \ h_1(T_1)>L+2l_\alpha,
\end{equation}
where $(u_1, h_1)$ is the solution of \eqref{p} with $u_0=\sigma_1\phi$, $v_\alpha$ is a solution of
\eqref{vala} and $l_\alpha$ is given in \eqref{la1}. The continuous dependence of the solution on
initial values yields that  there is a small $\epsilon>0$ such that the solution $(u_\epsilon, h_\epsilon)$
of \eqref{p} with $u_0=(\sigma_1-\epsilon)\phi$ satisfies \eqref{u-v_Z}. It then follows from Lemma
\ref{lem:spreadc} that spreading happens for $(u_\epsilon, h_\epsilon)$, which infers that
$\sigma_1-\epsilon\in\Sigma_1$. The comparison principle implies that $\sigma\in \Sigma_1$
for any $\sigma>\sigma_1-\epsilon$. Thus $\Sigma_1$ is an open interval.
\end{proof}

\begin{remark}\label{re:11}\rm
Denote $\sigma^*:=\inf\Sigma_1$. In view of the above lemma, we obtain
$$ \Sigma_1=(\sigma^*,\infty)\ \ \mbox{ and }\ \ \ \sigma^*\in[0,\infty).$$
\end{remark}

\subsection{Sufficient conditions for vanishing}
In this subsection we give some sufficient conditions for vanishing.
Firstly, let us consider the following eigenvalue problem:
 \begin{equation}\label{eigen-p}
 \left\{
 \begin{array}{ll}
 - \varphi'' -f'(0)\varphi =\lambda \varphi, & 0<x<L,\\
 - \varphi'' -g'(0)\varphi =\lambda \varphi, & x>L,\\
 \varphi'(0) = \varphi(\infty)=0, &  \\
 \varphi(L-0) = \varphi( L+0),\\
 \varphi'(L-0) = \varphi'( L+0).
 \end{array}
 \right.
 \end{equation}
Denote by $\lambda_1 (L)$ the principal eigenvalue of \eqref{eigen-p}. The existence and uniqueness of $\lambda_1 (L)$ is well-known; see \cite{DPS}.

Let us also consider the following eigenvalue problem on $\R$ associated with \eqref{eigen-p}:
 \begin{equation}\label{eigen-p1}
 \left\{
 \begin{array}{ll}
 - \varphi'' +h(x)\varphi =\lambda \varphi, & x\in\R,\\
 \varphi'(0) = \varphi(\pm \infty)=0,
 \end{array}
 \right.
 \end{equation}
where
 $$
h(x) = \left\{
\begin{array}{ll}
     -f'(0),\ \ & |x|\leq L,\\
      -g'(0), \ \ & |x|>L.
\end{array}
\right.
$$

As $h\in L^\infty(\R)$ and is symmetric with respect to the origin,
it is well known that the principal eigenvalue (or the so-called first eigenvalue) of \eqref{eigen-p1} exists and
coincides with that of problem \eqref{eigen-p}. Thus, we use $\lambda_1 (L)$ to denote the principal eigenvalue of
\eqref{eigen-p} and \eqref{eigen-p1}. The corresponding eigenfunction $\varphi_1^L$ of \eqref{eigen-p} satisfies
$\varphi_1^L\in C^1([0,\infty))\cap C^2([0,\infty)\setminus\{L\})$, $\varphi_1^L>0$ on $[0,\infty)$ and $(\varphi_1^L)'(0)=0$.

We recall the following assertion.

\begin{lem}[Lemma 3.1 of \cite{DPS}]\label{lem:1eigenvalue}
Let $\lambda_1 (L)$ be the principal eigenvalue of \eqref{eigen-p}. Then
\[
\lambda_1 (L)\in(-f'(0),-g'(0)),\ \ \mbox{ for any}\ L>0,
\]
and
\[
L=\frac{1}{\sqrt{f'(0)+\lambda_1(L)}}\arctan\sqrt{-\frac{g'(0)+\lambda_1(L)}{f'(0)+\lambda_1(L)}}.
\]
Moreover, $\lambda_1(L)$ is decreasing with respect to $L>0$, and $\lambda_1(L)<0$ if $L>L_*$, $\lambda_1(L)=0$ if
$L=L_*$, and $\lambda_1(L)>0$ if $0<L<L_*$, where $L_*=\frac{1}{\sqrt{f'(0)}}\arctan\sqrt{-\frac{g'(0)}{f'(0)}}$.
\end{lem}

Let $\lambda_1^R(L)$ be the principal eigenvalue of
\begin{equation}\label{eigen-p2}
\left\{
\begin{array}{ll}
- \varphi'' +h(x)\varphi =\lambda \varphi, & -R<x<R,\\
\varphi'(0) = \varphi(\pm R)=0.
\end{array}
\right.
\end{equation}
It follows from  \cite{BHR,BR} that
\begin{equation}\label{eig-1}
\mbox{$\lambda_1^R(L)$\ is decreasing in\ $R>0$}\ \, \mbox{and}\ \,\lim_{R\to\infty}\lambda_1^R(L)=\lambda_1 (L).
\end{equation}

In addition, some further properties of $\lambda_1^R(L)$ are collected as follows.

\begin{lem}\label{lem:1ei2} Let $L_*$ and $\lambda_1^R(L)$ be given as before.  We have
\begin{itemize}
\item[(i)] When $0<L\leq L_*$, then  $\lambda_1^R(L)>0$ for all $R>0$.

\item[(ii)] When $L> L_*$, then there is a unique positive constant $R^*:=R^*(L)$ such that  $\lambda_1^R(L)$ is negative $($resp. $0$,
or positive$)$ when $R>R^*$ $($resp. $R=R^*$, or $R<R^*$$)$. Moreover, $R^*$ is continuous and decreasing with respect to $L$.

\item[(iii)] Let $L_{**}:=\frac{\pi}{2\sqrt{f'(0)}}$. Then  $R^*(L)>L$ $($resp. $R^*(L)=L$, or $R^*(L)<L$$)$ when $L<L_{**}$
$($resp. $L=L_{**}$, or $L>L_{**}$$)$.
\end{itemize}
\end{lem}
\begin{proof}
(i) If $0<L\leq L_*$, it is known from Lemma \ref{lem:1eigenvalue} that $\lambda_1 (L)\geq 0$. By \eqref{eig-1}, we see that $\lambda_1^R(L)>0$ for all $R>0$.

(ii) If $L>L_*$, Lemma \ref{lem:1eigenvalue} implies that $\lambda_1 (L)< 0$. Due to \eqref{eig-1} and the continuity, there is a unique positive constant
$R^*$ such that  $\lambda_1^R(L)$ is negative (resp. $0$, or positive) when $R>R^*$ (resp. $R=R^*$, or $R<R^*$).
Moreover, $R^*$ is decreasing with respect to $L$. It's obvious that $R^*$ is continuous with respect to $L$.

(iii) Let us consider the continuous function
\[
\mathcal{J}(L):=R^*(L)-L.
\]
By the above discussion it follows that $\mathcal{J}(L)$ is decreasing and continuous with respect to $L>L_*$.
If $0<L\leq L_*$, then $\lambda_1^R(L)>0$ for all $R>0$, and so $R^*(L)\to +\infty$ as $L\to L_*+0$. It can be further checked that
$\lambda_1^R(L)=0$ if $L\geq \frac{\pi}{2\sqrt{f'(0)}}=R$, and in turn $R^*(L)=\frac{\pi}{2\sqrt{f'(0)}}$ for $L\geq \frac{\pi}{2\sqrt{f'(0)}}$.
As a result, we have
\[
\lim_{L\to L_*+0}\mathcal{J}(L)=+\infty,\ \ \ \lim_{L\to +\infty}\mathcal{J}(L)=-\infty.
\]
Thus there exists a unique $L_{**}\in(L_*,\infty)$ such that $\mathcal{J}(L_{**})=0$. Since $\mathcal{J}\big(\frac{\pi}{2\sqrt{f'(0)}}\big)=0$, hence
$L_{**}=\frac{\pi}{2\sqrt{f'(0)}}$. Thanks to the monotonicity of $\mathcal{J}(L)$ in $L$, we can derive all the other assertions of the lemma.
\end{proof}

\begin{remark}\label{rem12}\rm
In view of Lemma \ref{lem:1ei2}, we can define the following number
\begin{equation}\label{def:hh8}
h^* :=
\left\{
\begin{array}{ll}
  +\infty, & \mbox{when } 0\leq L\leq L_*,\\
  R^*(L), & \mbox{when } L>L_*.
\end{array}
\right.
\end{equation}
The number $h^*$ will paly an important role in our argument below.
\end{remark}

Next, based on the above results, we give the following sufficient condition for vanishing.

\begin{lem}\label{vfsma}
Let $L_{**}$ and $h^*$ be given in Lemma \ref{lem:1ei2} and \eqref{def:hh8}, respectively.
Assume that $0<L < L_{**}$, $L<h_0<h^*$ and $\|u_0\|_{L^\infty([0,h_0])}$ is
sufficiently small, then vanishing happens.
\end{lem}

\begin{proof}
For any given $h_1\in (h_0,h^*)$, we consider problem \eqref{eigen-p2} with $R =h_1$. Denote by
$\lambda^{h_1}_1$ and $\varphi$ with $\|\varphi\|_{L^\infty}=1$ the principal eigenvalue and the
corresponding positive eigenfunction of such a problem. Then $\lambda_1^{h_1}>0$ by Lemma \ref{lem:1ei2}.

Due to $\varphi(h_1)=0$, then $\varphi_x(h_1)\leq0$. Let $\zeta$ be denoted by the rightmost local
maximum point of $\varphi$ on $[0,h_1]$. Clearly, $\zeta\geq0$ and $\varphi_x(x)\leq0$ on $[\zeta,h_1]$. Set
\[
\delta: =\min\Big\{\frac{\lambda_1^{h_1}}{2},\ \ \frac{h_1}{h_0}-1, \ \ 1\Big\}, \ \
\eta:=\max\Big\{\zeta, \ h_1-\frac{\delta}{2}h_0 \Big\},\ \
\varepsilon_0:=\varphi(\eta),
\]
then $\varepsilon_0\leq 1$, and there exists $\varepsilon_1=\varepsilon_1(\delta)>0$
small such that
\[
-2\mu\varepsilon_1\varphi'(h_1)<\delta^2h_0,\ \ \ f(s)\leq (f'(0)+\delta)s,\ \ \ g(s)\leq (g'(0)+\delta)s\ \mbox{ for } s\in[0,\varepsilon_1].
\]
Define
\[
w(t,x) :=\varepsilon_0\varepsilon_1e^{-\delta t}\varphi(x)
\ \mbox{ for } (t,x)\in[0,\infty)\times[0,h_1].
\]
A direct calculation shows that for $(t,x)\in[0,\infty)\times[0,L)$,
\[
w_t-w_{xx}-f(w)\geq(\lambda_1^{h_1} - 2\delta)w\geq 0,
\]
and for $(t,x)\in[0,\infty)\times(L,h_1]$,
\[
w_t-w_{xx}-g(w)\geq(\lambda_1^{h_1} - 2\delta)w\geq 0.
\]
Since $\varphi(x)$ is the principal eigenfunction, clearly $w_x(t,0)=0$,
$w(t,L-0) = w(t, L+0)$ and  $w_x(t,L-0) = w_x(t, L+0),\,\forall t>0$. If we choose $u_0$ satisfying
\[
u_0(x) \leq\varepsilon_0 \varepsilon_1 \varphi(x) =  w(0,x)\ \ \mbox{ for } x\in [0,h_0],
\]
then it follows from Lemma \ref{lem:comp1}  that
\begin{equation*}\label{uw10}
u(t,x)\leq w(t,x)\ \mbox{ for }(t,x)\in[0,\tau )\times[0,h(t)],
\end{equation*}
where $\tau:=\sup\{t>0  : h(t)<h_1\}$. We shall prove that $\tau =\infty$. Once this is proved
we have $h(t)\leq h_1$ for all $t\geq 0$, and hence the vanishing conclusion follows.

To prove $\tau =\infty$, let us suppose for contradiction that $\tau<\infty$. Then $h(\tau)=h_1$. Inspired by \cite{SLZ}, we define
\[
\xi(t) :=h_0\Big(1+\delta-\frac{\delta}{2}e^{-\delta t}\Big),\
v(t,x) :=\varepsilon_1e^{-\delta t}\varphi(x-\xi(t)+h_1)
\]
for $t\geq 0,\ x\in I(t):= [\eta+\xi(t)-h_1, \ \xi(t)]$.

From the choice of $\eta$ it follows that for $t\geq0$,
$\eta+\xi(t)-h_1\geq \xi(t)-\frac{\delta}{2}h_0\geq h_0>L$.
By direct calculation, we see that for $t>0$ and $x\in I(t)$,
\[
v_t-v_{xx}-g(v)\geq (\lambda_1^{h_1}-2 \delta)v - \varepsilon_1e^{-\delta t}\xi'(t)\varphi' (x-\xi(t)+h_1)\geq  0,
\]
since $\xi'>0$ and $\varphi'(x-\xi(t)+h_1)<0$ for $t\geq 0$ and $x\in I(t)$. On the other hand,
by the choice of $\varepsilon_1$ we have
$$
\xi'(t)=\frac{\delta^2 h_0}{2} e^{-\delta t} \geq
-\mu\varepsilon_1 e^{-\delta t} \varphi'(h_1)
=- \mu v_x(t, \xi(t)).
$$

We claim that $h(t)\leq \xi(t)$ for all $t\in [0,\tau]$.  When $h(t)\leq \eta+ \xi(t)-h_1$
the claim is true since $\eta+ \xi(t)-h_1 <\xi(t)$.
Assume that the set $\{ 0\leq t\leq \tau  : h(t)>\eta+ \xi(t)-h_1\}\not= \emptyset$ consists of
some intervals and $[\tau_1, \tau_2]$ is one of them. Then $h(\tau_1 )=\eta+ \xi(\tau_1)-h_1$, and
on the left boundary $x=\eta+ \xi(t)-h_1$ of the domain $\Omega := \{(t,x) : t\in [\tau_1, \tau_2],\
\eta+ \xi(t)-h_1 \leq x \leq h(t)\}$ we have that for $t\in [\tau_1, \tau_2]$,
\begin{eqnarray*}
& & u(t, \eta+ \xi(t)-h_1)\leq  w(t, \eta+ \xi(t)-h_1)  =  \varepsilon_0 \varepsilon_1
e^{-\delta t} \varphi(\eta+ \xi(t)-h_1) \\
 & \leq &  \varepsilon_0 \varepsilon_1 e^{-\delta t} =\varepsilon_1 e^{-\delta t} \varphi(\eta)
  \equiv  v(t, \eta+ \xi(t)-h_1).
\end{eqnarray*}
Hence $v$ is an upper solution in $\Omega$ and by Lemma \ref{lem:comp1}, we have that $u\leq v$ in $\Omega$
and $h(t)< \xi(t)$ for $t\in [\tau_1, \tau_2]$.

In summary, our claim is proved and so
\[
h(\tau)\leq \xi (\tau)< \xi (\infty)\leq h_1,
\]
contradicting our assumption $h(\tau)=h_1$. This proves $\tau=\infty$, which ends the proof.
\end{proof}

Let $\phi \in \mathscr{X}(h_0)$ with $h_0 >L$, and $u_0=\sigma\phi$ with $\sigma>0$. Define
 $$
 \Sigma_0=\{\sigma_0: \mbox{ vanishing happens for } \sigma\in (0,\sigma_0]\}.$$
We end this subsection with some useful properties of the set $\Sigma_0$ for vanishing.
\begin{lem}\label{lem:sd0} Let $\Sigma_0$ be defined as before and $\sigma^*$ be given in Remark \ref{re:11}. The following assertions hold.
 \begin{itemize}
\item[\vspace{10pt}(i)] \;If $L> L_{**}$ or if $L_*<L<L_{**}$ and $h_0\geq R^*(L)$, then $\Sigma_0$ is empty;
\item[\vspace{10pt}(ii)] \;If $0<L\leq L_{*}$, then $\Sigma_0=(0,\sigma_*)$
                      for some $\sigma_*\in(0,\sigma^*)$;
\item[\vspace{10pt}(iii)] \;If $L_*<L< L_{**}$ and $h_0<R^*(L)$, then $\Sigma_0=(0,\sigma_*]$
                      for some $\sigma_*\in(0,\sigma^*)$.
\end{itemize}
\end{lem}

\begin{proof} (i)
Set $R:=h(1)$, then $\lambda_1^R(L)<0$ if $L>L_{**}$ or if $L_*<L<L_{**}$ and $h_0\geq R^*(L)$ by the properties of
$\lambda_1^R(L)$ in Lemma 3.2. The positive eigenfunction corresponding to $\lambda_1^R(L)$, denoted by $\varphi_1^{L,R}$,
solves \eqref{eigen-p2} and can be normalized so that $\|\varphi_1^{L,R}\|_{L^{\infty}}=1$. Let
 $$
\underline u(x)=\epsilon\varphi_1^{L,R}(x),\ \ x\in [0, R],
$$
where the constant $\epsilon>0$ can be chosen to be sufficiently small such that
\begin{equation}\nonumber
f(s)\geq\Big(f'(0)+\frac{\lambda_1^R(L)}{2}\Big)s\ \mbox{ and }\ g(s)\geq\Big(g'(0)+\frac{\lambda_1^R(L)}{2}\Big)s
\ \ \mbox{ for } s\in[0,\epsilon].
\end{equation}
As a result, for $t>2$ and $0\leq x \leq L$, it is easily seen that
$$
\underline u_t-\underline u_{xx}-f(\underline u) \leq \frac{\lambda_1^R(L)}{2}\underline u\leq 0,
$$
and for $t>2$ and $L<x\leq R$,
$$
\underline u_t-\underline u_{xx}-g(\underline u)\leq \frac{\lambda_1^R(L)}{2}\underline u\leq 0.
$$
One can check that $\underline u'(0)=0$,
$\underline u(L-0) = \underline u(L+0)$ and  $\underline u'(L-0) = \underline u'(L+0)$.

Moreover, since $u(2,x)>0$ for all $x\in[0, R]$, we can take $\epsilon$ to be smaller if necessary such that
$u(2,x)>\underline u(x)$ for all $x\in[0, R]$. Hence, $\underline u$ is a lower solution of \eqref{p}
for $t\geq 2,\,x\in[0, R]$. By the comparison principle, we obtain $u(t,x)\geq \underline u(x)$ for $t>2$ and $x\in[0, R]$.
This implies that vanishing can not happen for this $\sigma$, and thus $\Sigma_0$ is empty.

\smallskip

(ii) \ Since $0<L\leq L_{*}$, it follows from Lemma \ref{vfsma} and the parabolic comparison principle that
vanishing happens for all small $\sigma>0$, thus $\Sigma_0$ is not empty in this case.

Next we want to show that $\Sigma_0$ is an open interval. Fix any $\sigma_0\in \Sigma_0$, then vanishing happens for
$\sigma=\sigma_0$, and so for any $0<\delta\ll 1$, there exists $T_0>0$ large such that the solution $(u,h)$ of
\eqref{p} with $u_0=\sigma_0\phi$ satisfies
\[
\|u(t,\cdot)\|_{L^\infty ([0,h(T_0)])}<\delta\ \ \ \mbox{ and }\ \ \  h_\infty<\infty.
\]
By the continuous dependence of the solution of \eqref{p} on its initial values, we
can conclude that if $\varepsilon>0$ is sufficiently small, the solution $(u_\varepsilon,h_\varepsilon)$
of \eqref{p} with $u_0=(\sigma_0+\varepsilon)\phi$ satisfies
\[
\|u_\varepsilon(t,\cdot)\|_{L^\infty ([0,h_\varepsilon(T_0)])}<\delta\ \ \ \mbox{ and }\ \ \  h_\varepsilon(T_0)<h(T_0)+1.
\]
This, together with  Lemma \ref{vfsma} and the arbitrariness of $\delta$, yields that vanishing happens for
$(u_\varepsilon,h_\varepsilon)$, which implies that $\sigma_0+\varepsilon\in\Sigma_0$. Moreover, by the comparison
principle, $\sigma\in \Sigma_0$ for any $\sigma<\sigma_0+\varepsilon$. Thus $\Sigma_0$ is an open interval,
and so $\Sigma_0=(0, \sigma_*)$ with $\sigma_*\in(0,\sigma^*)$.

\smallskip

(iii) Since $L_*<L< L_{**}$ and $h_0<R^*(L)$, by Lemma \ref{vfsma} and the parabolic comparison principle we see
that vanishing happens for all small $\sigma>0$, thus $\Sigma_0$ is not empty. The definition of $\sigma_*$ implies
that $(0, \sigma_*)\subset\Sigma_0$. Then we intend to show that $\sigma_*\in\Sigma_0$. By Theorem \ref{thm:convergence},
it is suffice to prove that the solution $(u_*,h_*)$ of \eqref{p} with $u_0=\sigma_*\phi$ satisfies
\begin{equation}\label{eq:htoR}
h_{*,\infty}=R^*(L).
\end{equation}
We employ an indirect argument by assuming that $h_{*,\infty}<R^*(L)$ or $h_{*,\infty}>R^*(L)$. For the first case,
by the definition of vanishing and the  continuous dependence of the solution of \eqref{p} on the initial values, for
any small $\delta_1>0$, there is $T_1>0$ large such that if $\epsilon_1>0$ is sufficiently small, then the
solution $(u_1, h_1)$ of \eqref{p} with $u_0=(\sigma_*+\epsilon_1)\phi$, satisfies
\[
h_1(T_1)<R^*(L)\ \ \mbox{ and }\ \ u_1(T_1,x)<\delta_1
\mbox{ in } [0,h_1(T_1)],
\]
which implies that vanishing happens when $\sigma=\sigma_*+\epsilon_1$. This contradicts the
definition of $\sigma_*$. For the second case, we can find $T_2>0$ such that $h_*(T_2)>R^*(L)$.
By the continuous dependence of the solution of \eqref{p} on its initial values, we find
$\epsilon_2>0$ sufficiently small such that the solution of \eqref{p} with $u_0=(\sigma_*-\epsilon_2)\phi$,
denoted by $(u_2, h_2)$, satisfies
\[
h_2(T_2)>R^*(L).
\]
By the proof of (i) above, this implies that vanishing does not happen to $(u_2, h_2)$, a
contradiction to the definition of $\sigma_*$.

We claim that only when $\sigma=\sigma_*$, the solution $(u_*,h_*)$ of \eqref{p} with $u_0=\sigma_*\phi$
satisfies \eqref{eq:htoR}. Otherwise there exists $\tilde{\sigma}_1 \neq\sigma_*$
such that \eqref{eq:htoR} holds for $\sigma=\tilde {\sigma}_1$. Set $\tilde{\sigma} _2:=\sigma_*$, then the solution
of \eqref{p} with $u_0=\tilde {\sigma}_i\phi$, denoted by $(\tilde{u}_i,\tilde {h}_i)$ ($i=1,2$), satisfies
\[
\tilde{h}_{i, \infty}=R^*(L).
\]
Without loss of generality, we may suppose $\tilde{\sigma}_1>\tilde{\sigma}_2$. The comparison
principle yields that
\[
\tilde{h}_1(1)> \tilde{h}_2(1)\ \ \mbox{ and }\ \ \tilde{u}_1(1,x)>\tilde{u}_2(1,x)\ \mbox{in}\ [0,\tilde{h}_2(1)].
\]
Set
\[
\varepsilon_0=\sup\{\varepsilon>0:\tilde{u}_1(1,x)>\tilde{u}_2(1,x-
\varepsilon)\ \mbox{in}\ [\varepsilon,\tilde{h}_2(1)
+\varepsilon]\subset(0,\tilde{h}_1(1))\}
\]
and define
\[
\tilde{u}(t,x)=\tilde{u}_2(t+1,x-\varepsilon_0)\ \ \mbox{ and }\ \
\tilde{h}(t)=\tilde{h}_2(t+1)+\varepsilon_0.
\]
Clearly, $(\tilde{u}, \tilde{h})$ is the unique solution
of \eqref{p} with $u_0=\tilde{u}_2(1, x-\varepsilon_0)$ and $\tilde{h}_\infty=
\tilde{h}_{2,\infty}+\varepsilon_0$. The definition of $\varepsilon_0$
and the comparison principle conclude that
\[
\tilde{h}_{1,\infty}\geq \tilde{h}_\infty=\tilde{h}_{2,\infty}+\varepsilon_0>R^*(L).
\]
This leads to a contradiction against the definition of $\tilde{\sigma}_1$. Thus
\eqref{eq:htoR} holds only when $\sigma=\sigma_*$.
As a consequence, $\Sigma_0=(0, \sigma_*]$, which ends the proof of this lemma.
\end{proof}

\subsection{Proof of Theorem \ref{thm:dybe}}\label{subsect:pro1}
Based on the results derived in the previous subsections, we are going to give

\noindent
 {\bf Proof of Theorem \ref{thm:dybe}:}\ It follows from Remark \ref{re:11} that
 \begin{equation}\label{eq:se1}
 \Sigma_1=(\sigma^*,\infty)\ \mbox{ with }\ \sigma^*=\inf \Sigma_1\in[0,\infty).
 \end{equation}

 (I) When $0<L\leq L_*$, in view of Lemma \ref{lem:sd0}, we have $\Sigma_0=(0,\sigma_*)$ with
 $\sigma_*=\sup \Sigma_0\in(0,\sigma^*]$.  This, together with \eqref{eq:se1}, Theorem \ref{thm:convergence}
 and Lemma \ref{lem:staso} allow us to assert that each solution $u(t,x;\sigma\phi)$ with
 $\sigma\in [\sigma_*, \sigma^*]$  is a transition one.

\smallskip

(II) Let us assume that $L_*<L<\max\{L^*,\ L_{**}\}$ and divide our proof into four subcases.

Subcase (1): $L_*<L<L^*$ and $L^*< L_{**}$. We divide the initial
value $h_0$ into two cases: $h_0\geq R^*(L)$ and $h_0<R^*(L)$. In the first case, Lemma \ref{lem:sd0} (i) implies
that $\Sigma_0$ is empty, which means that vanishing does not happen for any
$\sigma>0$. This, combined with \eqref{eq:se1}, Theorem \ref{thm:convergence} and $L<L^*$, shows that each solution
$u(t,x;\sigma\phi)$ with $\sigma\in (0, \sigma^*]$ is a transition one. In the latter case, it follows from
Lemma \ref{lem:sd0} (iii) that $\Sigma_0=(0,\sigma_*]$ with $\sigma_*:=\sup \Sigma_0\in(0,\sigma^*]$.
On one hand, if $\sigma_*<\sigma^*$, by Theorem \ref{thm:convergence} and
Lemma \ref{lem:staso}, each solution $u(t,x;\sigma\phi)$ with $\sigma\in(\sigma_*,\sigma^*]$
is a transition one. On the other hand, if $\sigma_*=\sigma^*$, then the transition case does not happen.

Subcase (2): $L^*<L<L_{**}$ and $L^*< L_{**}$. The definition of $L^*$ implies that each solution
$u(t,x;\sigma\phi)$ with $\sigma>0$ is not a transition one.
There are two different cases: $h_0\geq R^*(L)$ and $h_0<R^*(L)$.
For the first case, by Lemma \ref{lem:sd0} (i), $\Sigma_0$ is empty, and in turn vanishing does not happen for any $\sigma>0$.
This, together with \eqref{eq:se1} and Theorem \ref{thm:convergence},
shows that spreading happens for all $\sigma>0$. For the latter case, we have from  Lemma \ref{lem:sd0}
(iii) that $\Sigma_0=(0,\sigma_*]$ with $\sigma_*=\sup \Sigma_0\in(0,\infty)$. Together with Theorem
\ref{thm:convergence} and \eqref{eq:se1}, we can conclude that $\sigma_*=\sigma^*$, and so $\Sigma_0=(0,\sigma^*]$ and
$\Sigma_1=(\sigma^*,\infty)$ in this case.

Subcase (3): $L_*<L<L_{**}$ and $L_{**}<L^*$. The proof is similar to that of subcase (1).

Subcase (4): $L_{**}<L<L^*$ and $L_{**}<L^*$. By Lemma \ref{lem:sd0} (i), we obtain that vanishing does not happen
for any $\sigma>0$. This, combined with Theorem \ref{thm:convergence} and \eqref{eq:se1},
shows that each solution $u(t,x;\sigma\phi)$ with $\sigma\in (0, \sigma^*]$ is a transition one if $\sigma^*>0$.

 \smallskip

(III) Let us consider the case $L>\max\{L^*,\ L_{**}\}$. As $L>L_{**}$, vanishing does not happen for problem \eqref{p}
due to Lemma \ref{lem:sd0} (i). The definition of $L^*$, together with Theorem \ref{thm:convergence} and Lemma \ref{lem:staso},
imply that only spreading can happen for problem \eqref{p} with $L>\max\{L^*,\ L_{**}\}$.

The proof of Theorem \ref{thm:dybe} is now complete. {\hfill $\Box$}

\subsection{Asymptotic profiles of spreading solutions} \label{subsect:asspeed}
In this subsection we study the asymptotic profiles of spreading solutions and prove Theorem \ref{thm:profile of spreading sol}.
Let us first state the following known result from \cite[Theorem 6.2]{DuLou},
which will play an important role in our analysis of the asymptotic spreading speed.

\begin{lem}\label{lem:semi-wave}
Let $c_0$ be the speed of traveling wave of $u_t=u_{xx}+g(u)$. For any $c\in[0,c_0)$, the following problem
\begin{equation}\label{psemi-wave}
\left\{
 \begin{array}{ll}
 q_c''-cq_c'+g(q_c) =0,\ \ q_c(z)>0, & 0< z< \infty,\\
 q_c(0)=0,\ \ q_c(\infty)=1,
  \end{array}
 \right.
\end{equation}
admits a unique positive solution $q_c(z)$ with $q_c'(z)>0$ for $z\geq0$. Moreover, for each $\mu>0$,
there exists a unique $c^*=c_{\mu}^*\in(0,c_0)$ such that \eqref{psemi-wave} has a unique solution
pair $(c, q)=(c^*, q_{c^*})$ satisfying $c^*=\mu q_{c^*}'(0)$,
and $c^*_\mu$ is increasing in $\mu$ with $c^*_\mu\to c_0$ as $\mu\to\infty$.
\end{lem}

Based on the above result, we will prove the boundedness of $h(t)-c^*t$ and that $u(t,\cdot) \approx 1$
in the domain $[0, h(t)-X]$, where $X>0$ is a large number to be determined.

\begin{prop}\label{pro:sigma01}
Assume that spreading happens for the solution $(u,h)$. Let $(c^*, q_{c^*})$ be the unique solution
pair of \eqref{psemi-wave}. Then
\begin{itemize}
\item[(i)] there exists $C>0$ such that
\begin{equation}\label{hh1}
|h(t)-c^*t |\leq  C \ \ \mbox{for all } t\geq 0 ;
\end{equation}

\item[(ii)] for any small $\epsilon>0$, then there exist $X_\epsilon >0$ and $T_\epsilon >0$ such that
\begin{equation}\label{u to P near 0}
\|u(t,\cdot ) - 1 \|_{L^\infty ([0, h(t) -X_\epsilon])} \leq  \epsilon \ \mbox{ for } t> T_\epsilon.
\end{equation}
\end{itemize}
\end{prop}

\begin{proof} For clarity we divide the proof into three steps.

\smallskip

{\bf Step 1}. Estimate the upper bounds for $h(t)$ and $u(t,x)$.

Since $(c^*, q_{c^*})$ is a solution of \eqref{psemi-wave} with $q_{c^*}(z)\to 1$
as $z\to\infty$, a simple analysis on the $q-q'$ phase plane around the point $(1,0)$ gives
\[
q'_{c^*}(z)=[-2\gamma+o(1)](q_{c^*}(z)-1),\ \ \mbox{as}\ z\to\infty,
\]
where
\[
\gamma:=\frac{1}{4}[\sqrt{(c^*)^2-4g'(1)}-c^*]>0.
\]
Then we can find $X_0>0$ large such that there exists $K_0>0$ satisfying
\begin{equation}\label{eq:qc1}
q_{c^*}(z)\geq 1-K_0 e^{-\gamma z}\ \ \mbox{for }  z\geq X_0.
\end{equation}
Choosing $0<\delta_1 <\frac{1}{2}\min\{\gamma c^*,-g'(1)\}$, it follows from  {\bf{(H)}} that there
is $\varepsilon>0$ such that
\begin{equation}\label{fuzhu3}
f'(v),\ g'(v) \leq -\delta_1 \    \mbox{ for } \ v\in[1-\varepsilon ,1+\varepsilon].
\end{equation}

Let $v(t)$ be the solution of $v_t=F(v)$
with initial value $v(0)=\|u_0\|_\infty +1$, where
\[
F(v)=\max\{f(v),\ g(v)\}.
\]
Due to $F(v)<0$ for $v>1$, the function $v(t)$ decreases to $1$ as $t\to \infty$. Hence, for
$\varepsilon>0$ in \eqref{fuzhu3}, there exists $t_0>0$ large such that $1<v(t)<1+\varepsilon$
for $t \geq t_0$. By \eqref{fuzhu3} we have $v_t = F(v) \leq \delta_1 (1-v)$ for $t\geq t_0$.
Clearly $v(t)$ is an upper solution of \eqref{p}, and so
\begin{equation}\label{Vtoto1}
u(t,x) \leq  v (t)\leq 1+ M_0  e^{-\delta_1 t}\ \ \ \mbox{ for }\
t \geq t_0,\ 0\leq x\leq h(t),
\end{equation}
with $M_0:= \varepsilon e^{\delta_1 t_0}$. Take $T' >t_0$ large such that
\begin{equation}\label{equ:ee1}
h(t)\geq 2L\ \mbox{ and }\ M_0 e^{-\delta_1t}< \varepsilon /2\ \mbox{ for }\ t\geq T'.
\end{equation}
Since $q_{c^*}(z)\to 1$ as $z\to\infty$, we can find $X >X_0$ large such that, with
$M'=2M_0$,
\begin{equation}\label{U1a1}
(1+M'e^{-\delta_1 T'  })q_{c^*}(z)\geq 1+ M_0 e^{-\delta_1 T'}
\ \mbox{ for }  \ z\geq X.
\end{equation}

Now we construct a finer upper solution $(\bar{u}, \bar{h} )$ to \eqref{p} as follows:
\begin{align*}
& \bar{h}(t): =c^*t+ h (T')+ K M'(e^{-\delta_1 T'}-e^{-\delta_1 t})+X+2L \ \
\mbox{ for } t\geq T',\\
& \bar{u}(t,x):=(1+M'e^{-\delta_1 t})q_{c^*}(\bar{h}(t)-x)\  \ \mbox{ for } t\geq T',\ 2L\leq x \leq \bar{h}(t),
\end{align*}
where $K$ is a positive constant to be determined below.
Clearly, for all $t\geq T'$, we have $\bar{u} (t, \bar{h}(t))=0$, and
\[
-\mu \bar{u} _x(t,\bar{h} (t))=\mu(1+M'e^{-\delta_1 t})q'_{c^*}(0)=(1+M'e^{-\delta_1 t})c^* <
c^*+M' K \delta_1 e^{-\delta_1 t} = \bar{h} '(t)
\]
provided that we choose $K$ with $K\delta_1 > c^*$.  By the definition of $\bar{h}$, we
have $h (T')<\bar{h}(T')$.
It follows from \eqref{eq:qc1}-\eqref{U1a1} that
for $x\in[2L, h (T')]$,
\begin{eqnarray*}
\bar{u} (T',x) &  = & (1+M'e^{-\delta_1 T'})q_{c^*}(c^*T'+h (T')+X+2L -x) \\
& \geq & 1+ M_0 e^{-\delta_1 T'} \geq  u(T',x),
\end{eqnarray*}
and for $t\geq T'$,
\begin{eqnarray*}
\bar{u} (t , 2L) &  \geq & (1+M'e^{-\delta_1 t})q_{c^*}(c^*t+X) \\
&\geq & (1+M'e^{-\delta_1 t})(1-K_0 e^{-\gamma (c^*t+X)})\\
& \geq & 1+ 2M_0 e^{ - \delta_1 t }-(1+\varepsilon)K_0 e^{-\gamma (c^*t+X)}\\
& \geq & 1+ M_0 e^{-\delta_1 t} \geq  u(t ,2L),
\end{eqnarray*}
by the fact that for $\delta_1\leq \frac{\gamma c^*}{2}$, due to $X\gg 1$ and $T'\gg 1$, then for $t\geq T'$,
\[
M_0 e^{ - \delta_1 t }-(1+\varepsilon)K_0 e^{-\gamma (c^*t+X)}\geq M_0 e^{ - \delta_1 t }[1-e^{(\delta_1- \gamma c^*)t}]\geq 0.
\]

We now show that
\begin{equation*}\label{u+ upper}
 \mathcal{N} \bar{u} := \bar{u}_t - \bar{u}_{xx} -g(\bar{u}) \geq 0\quad \mbox{ for }\ t>T',\ x\in [2L, \bar{h}(t)].
\end{equation*}
In fact, by some direct calculation we find that
\begin{eqnarray*}
\mathcal{N} \bar{u}  & = & M'e^{-\delta_1 t}\big\{g(q_{c^*})+K \delta_1 (1+M'e^{-\delta_1 t})q'_{c^*}
-\delta_1 q_{c^*}\big\} + g(q_{c^*})- g((1+M'e^{-\delta_1 t})q_{c^*})\\
 & = & \mathcal{F}_1:= M'e^{-\delta_1 t}\big\{g(q_{c^*})+K \delta_1 (1+M'e^{-\delta_1 t})q'_{c^*}-[g'((1+ \rho_1 M'e^{-\delta_1 t})q_{c^*})+\delta_1]q_{c^*}\big\}
\end{eqnarray*}
for some $\rho_1 \in (0,1)$.  Since $q_{c^*}(z)\to 1$ as $z\to \infty$, there is $z_0>0$ such that
\[
q_{c^*}(z)\geq 1-\varepsilon \ \ \mbox{ for }\ z\geq z_0.
\]
When $t> T' $ and $\bar{h} (t)-x>z_0$, clearly $\mathcal{F}_1\geq 0$ by \eqref{fuzhu3} and the fact that
$M'e^{-\delta_1 t}\leq\varepsilon $ for $t> T' $.
When  $t> T' $ and  $0\leq \bar{h}(t)-x\leq z_0$, we have
$$
\mathcal{F}_1  \geq  M'e^{-\delta_1 t}(K \delta_1 D_1 - D_2 -\delta_1-D_3)\geq 0,\quad \mbox{ provided } K>0 \mbox{ is sufficiently large},
$$
where
$$
D_1:=\min_{z\in[0,z_0]}q'_{c^*}(z)>0, \ \ D_2:=\max_{s\in[0,1+\varepsilon]}g'(s)\ \ \mbox{ and } \ D_3:=\max_{s\in[0,1]}g(s).
$$

In summary, $(\bar{u}, \bar{h})$ is an upper solution of \eqref{p}.  It follows from Lemma \ref{lem:comp2} that
$$
h(t) \leq \bar{h}(t) \mbox{ for } t>T'\quad \mbox{and } \ \  u(t,x)\leq \bar{u}(t,x) \leq 1+M' e^{-\delta t},\quad  t>T',\  x\in[2L, h(t)].
$$
By the definition of $\bar{h}$, for $C_r := h(T')+X +KM'+2L+\max\limits_{t\in[0,T']}|\overline{h}(t)-h(t)|$, we see that
\begin{equation}\label{hbd}
h(t)< c^*t +C_r \ \mbox{ for all } t\geq 0.
\end{equation}
For any $\epsilon>0$, by choosing $T_1(\epsilon) >T'$ large such that
$M' e^{-\delta T_1(\epsilon)} < \epsilon$, thanks to the definition of $\bar{u}$, we have
\begin{equation}\label{v<1+epsilon/P}
u(t,x)\leq \bar{u}(t,x) \leq 1 +  \epsilon\quad \ \mbox{ for }\   t> T_1(\epsilon),\ x\in [2L, h(t)].
\end{equation}

\smallskip

{\bf Step 2}. Estimate the lower bounds for $h(t)$ and $u(t,x)$.

It follows from  \cite[Lemma 6.5]{DuLou} that for any given $c\in (0,c^*)$
there exist $\delta_2\in(0,-g'(1))$, $T_1>0$
and $ M>0 $  such that for $t\geq T_1$,
\begin{equation}\label{eq:cah}
 h(t)\geq ct\ \ \mbox{ and }\ \  u(t,x)\geq 1-Me^{-\delta_2 t}\quad \mbox{for } x\in [0, ct].
\end{equation}
Since $\delta_2\in(0,-g'(1))$,  we can find some $\eta>0$ small such that
\begin{equation}\label{eq:gg1}
g'(v)\leq -\delta_2\ \ \ \mbox{ for }\ v\in[1-\eta,1+\eta].
\end{equation}
Moreover, we define constants $z_\eta$ and $C_\eta$ as follows:
\[
q_{c^*}(z_\eta)=1-\frac{\eta}{2},\ \ \ \ C_\eta=\min_{0\leq z\leq z_\eta} q'_{c^*}(z)>0.
\]
Then we take $T''>T_1$ such that
\begin{equation}\label{eq:me1}
ct\geq 2L\ \mbox{ and }\ Me^{-\delta_2 t}\leq \frac{\eta}{2}\ \ \mbox{ for }\ t\geq T''.
\end{equation}
Define
\begin{align*}
& \underline{h}(t): =c^*(t-T'')+ cT''- K_1 M(e^{-\delta_2 T''}-e^{-\delta_2 t}) \ \
\mbox{ for } t\geq T'',\\
& \underline{u}(t,x):=(1-Me^{-\delta_2 t})q_{c^*}(\underline{h}(t)-x)\  \ \mbox{ for } t\geq T'' ,\ 2L\leq x \leq \underline{h}(t),
\end{align*}
where $K_1$ is a positive constant to be determined below.
Clearly, for all $t\geq T'' $, we infer
\[
\underline{u} (t, \underline{h}(t))=0,\ \ \mbox{ and }\ \ \ \underline{u} (t, 2L)\leq 1-Me^{-\delta_2 t}\leq u(t, 2L),
\]
\[
-\mu \underline{u} _x(t,\underline{h} (t))=\mu(1-Me^{-\delta_2 t})q_{c^*}(0)=(1-Me^{-\delta_2 t})c^* >c^*-K_1M \delta_2 e^{-\delta_2 t} = \underline{h} '(t)
\]
provided that we choose $K_1$ with $K_1\delta_2 > c^*$.  By the definition of $\underline{h}$, we have $h (T'')>\underline{h}(T'')$.
In view of \eqref{eq:cah}, we further have
for $x\in[ 2L, \underline{h }(T'')]$,
\begin{eqnarray*}
\underline{u} (T'',x) \leq 1-Me^{-\delta_2 T''}\leq u(T'',x).
\end{eqnarray*}

We now show that
\begin{equation*}\label{u+ upper}
 \mathcal{N} \underline{u} := \underline{u}_t - \underline{u}_{xx} -g(\underline{u}) \leq 0
 \quad \ \mbox{ for }\ \ t>T'',\ x\in [2L, \underline{h}(t)].
\end{equation*}
In fact, direct calculation gives
\begin{eqnarray*}
\mathcal{N} \underline{u}
 & = & -Me^{-\delta_2 t}\big\{g(q_{c^*})+K_1 \delta_2 (1-Me^{-\delta_2 t})q'_{c^*} -\delta_2 q_{c^*}\big\} + g(q_{c^*})- g((1-Me^{-\delta_2 t})q_{c^*})\\
 & = & \mathcal{F}_2:=- Me^{-\delta_2 t}\big\{g(q_{c^*})+K_1 \delta_2 (1-Me^{-\delta_2 t})q'_{c^*}-[g'((1-\rho_2 Me^{-\delta_2 t})q_{c^*})+\delta_2]q_{c^*}\big\},
\end{eqnarray*}
for some $\rho_2 \in (0,1)$.

When  $t> T''$ and $\underline{h} (t)-x>z_\eta$, $\mathcal{F}_2 \leq 0$ by the fact that for $z\geq z_\eta$, due to \eqref{eq:me1},
\[
1-\eta \leq \big[1-e^{-\delta_2 t}\big]q_{c^*}(z)\leq \big[1-\rho_2 e^{-\delta_2 t}\big]q_{c^*}(z)\leq 1,
\]
and hence, by \eqref{eq:gg1}, $g'((1-\rho_2 Me^{-\delta_2 t})q_{c^*})+\delta_2\leq 0$.

When  $t> T''$ and $0\leq \underline{h} (t)-x\leq z_\eta$, we have
\begin{eqnarray*}
\mathcal{F}_2 & \leq & - Me^{-\delta_2 t}\Big\{\min_{0\leq s\leq 1}g(s)+K_1 \delta_2 (1-Me^{-\delta_2 t})q'_{c^*}-\max_{0\leq s\leq 1}g'(s)-\delta_2\Big\}\\
& \leq & Me^{-\delta_2 t}(1-Me^{-\delta_2 t})\Big\{\frac{\max_{0\leq s\leq 1}g'(s)+\delta_2-\min_{0\leq s\leq 1}g(s)}{1-Me^{-\delta_2 t}}-K_1\delta_2 C_\eta\Big\}\\
& \leq & Me^{-\delta_2 t}(1-Me^{-\delta_2 t})\Big\{\frac{\max_{0\leq s\leq 1}g'(s)+\delta_2-\min_{0\leq s\leq 1}g(s)}{1-Me^{-\delta_2 T''}}-K_1\delta_2 C_\eta\Big\}.
\end{eqnarray*}
By taking $K_1>0$ sufficiently large, we have $\mathcal{F}_2\leq 0$.

Consequently, $(\underline{u}, \underline{h})$ is a lower solution. It follows from the comparison principle that
\begin{equation*}\label{h>h->R}
h(t) \geq \underline{h}(t)  \mbox{ for } t >T'',
\quad u(t,x)\geq  \underline{u}(t,x) \mbox{ for }\ t>T'',\ x\in [2L, \underline{h}(t)].
\end{equation*}
Hence
\begin{equation}\label{hbd2}
h(t)\geq \underline{h}(t) - \max_{t\in[0,T'']}|h(t)-\underline{h}(t)| \geq
c^*t -C_l \ \ \mbox{ for all } t\geq 0,
\end{equation}
where $C_l = \underset{t\in[0,T'']}{\max}|h(t)-\underline{h}(t)|+c^*T'' + K_1 M$.
This, together with \eqref{hbd}, verifies \eqref{hh1}.

Moreover, for any $\epsilon>0$, since $q_{c^*}(z) \to 1$ as $z\to \infty$,
there exists $X_1(\epsilon)>0$ such that
$$
q_{c^*}(z)> 1- \epsilon/2\ \ \ \mbox{  for }\ \  z\geq X_1(\epsilon).
$$
For $(t,x)\in \Omega_1 := \{ (t,x) : t>T'',\ 2L\leq x\leq h(t) -C_r -C_l -X_1(\epsilon)\}$, one can get from
\eqref{hbd} and \eqref{hbd2} that
$$
\underline{h}(t) -x \geq c^*t -C_l -x \geq h(t) - C_r -C_l  -x \geq X_1(\epsilon),
$$
and hence
$$
u(t,x) \geq \underline{u}(t,x) \geq (1-M e^{-\delta_2 t} ) q_{c^*}(X_1(\epsilon)) \geq
(1-M e^{-\delta_2 t} )( 1- \epsilon/2) \mbox{ for } (t,x)\in \Omega_1.
$$

Let us choose $T_2(\epsilon) >T''$ such that $2 M e^{-\delta_2 T_2(\epsilon)} <\epsilon$, then
\begin{equation}\label{v>1-epsilonP}
u(t,x)\geq ( 1- \epsilon/2)^2 > 1-\epsilon\ \ \ \mbox{ for }
(t,x)\in \Omega_1 \mbox{ and } t>T_2(\epsilon).
\end{equation}

\smallskip

{\bf Step 3}. Complete the proof of \eqref{u to P near 0}.

Since spreading happens for $(u,h)$, then
there exists $T_0>0$ large such that
\begin{equation}\label{eq:utt1}
|u(t,x)-1| \leq \epsilon \ \ \  \mbox{ for }  t>T_0,\ 0\leq x\leq 2L.
\end{equation}
Denote $T_\epsilon := \max\{T_0,\ T_1(\epsilon),\ T_2(\epsilon)\}$ and $X_\epsilon := C_r + C_l +X_1(\epsilon)$. Then
by \eqref{v<1+epsilon/P}, \eqref{v>1-epsilonP} and \eqref{eq:utt1},  we have
$$
|u(t,x)-1| \leq \epsilon \ \ \ \mbox{ for } t>T_\epsilon,\ 0\leq x\leq h(t) -X_\epsilon.
$$
This yields the estimate in \eqref{u to P near 0}, which ends the proof.
\end{proof}

Making use of the above results that have already been proved, we are going to give

\smallskip

\noindent
 {\bf Proof of Theorem \ref{thm:profile of spreading sol}:}\
We first prove \eqref{HWt1}. It follows from \eqref{hh1} in Proposition \ref{pro:sigma01} that there is $C>0$ such that
\[
-C\leq h(t)-c^*t\leq C \ \ \mbox{ for }\ t\geq 0.
\]
We define
\[
\tilde{u}(t,y):=u(t,y+c^*t-2C),\ \ \ \tilde{h}(t):=h(t)-c^*t+2C, \ \ t\geq0.
\]
Let $t_n\to\infty$ be an arbitrary sequence and define
\[
\tilde{u}_n(t,y):=\tilde{u}(t+t_n,y),\ \ \ \tilde{h}_n(t):=\tilde{h}(t+t_n).
\]
By a similar argument as that in \cite[Lemmas 4.2 and 4.6]{DMZ2}, by passing to a subsequence
if necessary, we obtain that there is a constant $H_0\in \R$ such that
\[
\tilde{h}_n\to H_0\ \mbox{ in }\ C^{1+\frac{\gamma}{2}}_{loc}(\R),
\]
with $\gamma\in(0,1)$. The arbitrariness of $\{t_n\}$ and the definition of $\tilde{h}_n$ imply that
\[
h'(t)\to c^*,\ \ \mbox{as}\ t\to\infty.
\]

\medskip
We next prove \eqref{WHt1}. To this end, we use the moving coordinate $z:= x-h(t)$. Set
$$
v(t,z) := u(t, z+h(t))\ \  \mbox{ for }\ z\in [2L-h(t), 0],\ t\geq 0.
$$
Then $v$ solves
\begin{equation*}\label{p u2}
\left\{
\begin{array}{ll}
 v_t =v_{zz}+h'(t)v_z+ g(v), &  2L-h(t)<z<0,\ t>0,\\
 v(t, 0)= 0,\ h'(t) = -\mu v_z(t,0), &  t>0,\\
 v(t, 2L-h(t))= u(t,2L), &  t>0,
 \end{array}
\right.
\end{equation*}
where, as $t\to\infty$, $2L-h(t)\to-\infty$ and $v(t,2L-h(t))\to 1$.

Consider the $\omega$-limit set $\omega(v)$ of $v(t,\cdot)$ in the topology of
$C^2_{loc}((-\infty,0])$. As $v$ is bounded in $L_{loc}^\infty$ norm, then $\omega(v)$
is not empty. It follows from the ideas developed by Du and Matano \cite{DM} and Du and Lou \cite{DuLou}
that $\omega(v)$ consists of only solutions of
$w_{zz}+c^*w_z+g(w)=0$ for $z\in(-\infty,0]$ in virtue of $h'(t)\to c^*$
as $t\to\infty$. For each $w\in \omega(v)$, we have
\[
w(-\infty)=1, \ \ \ w(0)=0, \ \ \mbox{ and }\ \ c^*=-\mu w_z(0).
\]
Thus $\omega(v)=\{q_{c^*}(-z)\}$, which in turn implies that, for any $C>0$,
$$
\|v(t, \cdot) -q_{c^*}(\cdot)\|_{L^\infty ([-C,\ 0])} \to 0, \quad \mbox{as}\ t\to \infty,
$$
or, equivalently,
\begin{equation}\label{u to U near h(t)}
\|u(t, \cdot) - q_{c^*}(h(t) -\cdot)\|_{L^\infty ([h(t)-C,\ h(t)])} \to 0, \quad \mbox{as}\ t\to \infty.
\end{equation}

For any given small $\epsilon >0$, it follows from \eqref{u to P near 0} in Proposition \ref{pro:sigma01} that there exist
two positive constants $X_\epsilon$ and $T_\epsilon$ such that
$$
|u(t,x) - 1| \leq \epsilon\quad \mbox{for }  t>T_\epsilon,\ 0\leq x\leq h(t) - X_\epsilon.
$$
Since $q_{c^*}(\infty)=1$, there exists $X^*_\epsilon > X_\epsilon$ such that
$$
|q_{c^*}(h(t) -x) -1| \leq \epsilon\quad\  \mbox{ for } x\leq h(t) -X^*_\epsilon.
$$
Combining the above two inequalities, we deduce
$$
|u(t,x) - q_{c^*}(h(t) -x) | \leq 2\epsilon\quad\ \mbox{ for }\ t> T_\epsilon,\ 0\leq x\leq h(t)-X^*_\epsilon.
$$
Taking $C= X^*_\epsilon$ in \eqref{u to U near h(t)} we see that for some $T^{*}_\epsilon >T_\epsilon$,
$$
|u(t, x) -q_{c^*}(h(t)-x)| \leq \epsilon\ \quad \mbox{ for }\ t>T^{*}_\epsilon, h(t) -X^*_\epsilon \leq x \leq h(t).
$$
This proves \eqref{WHt1}. Thus, Theorem \ref{thm:profile of spreading sol} is verified. {\hfill $\Box$}

\section{A Separated Protection Zone: proof of Theorem \ref{thm:dybe1}}\label{sect:multiple case}
In this section, we consider system \eqref{q} and prove Theorem \ref{thm:dybe1} in the same
spirit as that of Theorems \ref{thm:dybe} and \ref{thm:profile of spreading sol}. Consider the following eigenvalue problem:
\begin{equation}\label{eigen-q}
\left\{
 \begin{array}{ll}
 - \varphi'' -f'(0)\varphi =\lambda \varphi, & x \in(L_1, L_2),\\
 - \varphi'' -g'(0)\varphi =\lambda \varphi, & x\in(0, L_1)\cup(L_2,\infty),\\
 \varphi'(0)=\varphi(\infty)=0, \\
 \varphi(L_i-0) = \varphi( L_i+0), & i=1, 2, \\
 \varphi'(L_i-0) = \varphi'( L_i+0), & i=1, 2.
 \end{array}
 \right.
\end{equation}

Let us also consider the following eigenvalue problem on $\R$ associated with \eqref{eigen-q}:
 \begin{equation}\label{eigen-q1}
 \left\{
 \begin{array}{ll}
 - \varphi'' +\tilde h(x)\varphi =\lambda \varphi, & x\in\R,\\
 \varphi'(0) = \varphi(\pm \infty)=0,
 \end{array}
 \right.
 \end{equation}
where
 $$
\tilde h(x) = \left\{
\begin{array}{ll}
     -f'(0),\ \ & |x|\in[L_1, L_2],\\
      -g'(0), \ \ & |x|\in[0, L_1)\cup (L_2, \infty).
\end{array}
\right.
$$

It is well known that \eqref{eigen-q} and \eqref{eigen-q1} have the same principal eigenvalue,
denoted by $\tilde \lambda_1 (L)$. From \cite{BHR,BR} we also observe that
 \begin{equation}\label{eq:rtoq}
 \mbox{$\tilde\lambda_1^R(L)$\ is decreasing in\ $R>0$}\ \, \mbox{and}\ \,\lim_{R\to\infty}\tilde\lambda_1^R(L)=\tilde\lambda_1 (L),
 \end{equation}
where $\tilde\lambda_1^R(L)$ is the principal eigenvalue of
 \begin{equation*}\label{eigen-q2}
 \left\{
 \begin{array}{ll}
 - \varphi'' +\tilde h(x)\varphi =\lambda \varphi, & -R<x<R,\\
 \varphi'(0) = \varphi(\pm R)=0.
 \end{array}
 \right.
 \end{equation*}

We have the following results.
\begin{lem}[Lemma 4.1 of \cite{DPS}]\label{lem:1eigenvaluemp2}
For any given $0<L_1<L_2$, let $L=L_2-L_1$ and $\tilde{\lambda}_1(L)$
be the principal eigenvalue of \eqref{eigen-q}. Then we have
\[
\tilde{\lambda}_1(L)\in(-f'(0),-g'(0)),
\]
and
\[
L=\frac{1}{\theta_2}\Big\{ \arctan\Big[\frac{\theta_1}{\theta_2}\cdot\frac{e^{\theta_1 L_1}
-e^{-\theta_1 L_1}}{e^{\theta_1 L_1}+e^{-\theta_1 L_1}}\Big]+  \arctan \frac{\theta_1}{\theta_2}\Big\},
\]
where $\theta_1=\sqrt{-(g'(0)+\tilde{\lambda}_1(L))}$ and $\theta_2=\sqrt{f'(0)+\tilde{\lambda}_1(L)}$.
Moreover, $\tilde{\lambda}_1(L)$ is decreasing with respect to $L>0$, and there is a unique
$\tilde{L}_*>L_*$ such that $\tilde{\lambda}_1(L)<0$ if $L>\tilde{L}_*$, $\tilde{\lambda}_1(L)=0$
if $L=\tilde{L}_*$, and $\tilde{\lambda}_1(L)>0$ if $0<L<\tilde{L}_*$.
\end{lem}

\begin{lem}\label{lem:s1ei2}
For any given $0<L_1<L_2$, set $L=L_2-L_1$. Let $\tilde{L}_*$ and $\tilde{\lambda}_1^R(L)$ be given as above.
The following assertions hold.
\begin{itemize}
\item[(i)] If $0<L\leq \tilde{L}_*$, then  $\tilde{\lambda}_1^R(L)>0$ for all $R>0$.

\item[(ii)] If $L> \tilde{L}_*$, there is a unique positive constant $\tilde{R}^*:=\tilde{R}^*(L)$ such that
$\tilde{\lambda}_1^R(L)$ is negative (resp. $0$, or positive) when $R>\tilde{R}^*$ (resp. $R=\tilde{R}^*$, or $R<\tilde{R}^*$).
For any given $L_1>0$, $\tilde{R}^*$ is continuous and decreasing with respect to $L$.

\item[(iii)] For any given $L_1>0$, let
\begin{equation}\label{ee1}
\tilde{L}_{**}:=\frac{1}{\sqrt{f'(0)}}\Big\{ \arctan\Big[\sqrt{-\frac{g'(0)}{f'(0)}}\cdot\frac{e^{\sqrt{-g'(0)} L_1}
-e^{-\sqrt{-g'(0)} L_1}}{e^{\sqrt{-g'(0)} L_1}+e^{-\sqrt{-g'(0)} L_1}}\Big]+  \frac{\pi}{2}\Big\}>L_{**},
\end{equation}
then $\tilde{R}^*(L)>L_2$ (resp. $\tilde{R}^*(L)=L_2$, or $\tilde{R}^*(L)<L_2$) when $L<\tilde{L}_{**}$
(resp. $L=\tilde{L}_{**}$, or $L>\tilde{L}_{**}$).
\end{itemize}
\end{lem}
\begin{proof}
(i) If $0<L\leq \tilde{L}_*$, by Lemma \ref{lem:1eigenvaluemp2}, $\tilde{\lambda}_1 (L)\geq 0$. Thanks to \eqref{eq:rtoq},
we see that $\tilde{\lambda}_1^R(L)>0$ for all $R>0$.

(ii) If $\tilde{L}_*<\ L$, by Lemma \ref{lem:1eigenvaluemp2}, $\tilde{\lambda}_1 (L)< 0$. This, together with \eqref{eq:rtoq},
implies that there is a unique constant $\tilde{R}^*>0$ such that  $\tilde{\lambda}_1^R(L)$ is negative (resp. $0$,
or positive) when $R>\tilde{R}^*$ (resp. $R=\tilde{R}^*$, or $R<\tilde{R}^*$). Moreover, one can check that $\tilde{R}^*$ is
continuous and decreasing with respect to $L$.

(iii) For any given $L_1>0$, let us consider the following function
\[
\tilde{\mathcal{J}}(L):=\tilde{R}^*(L)-L_1-L.
\]
By the above discussion we know that for any given $L_1>0$, $\tilde{\mathcal{J}}(L)$ is decreasing and
continuous with respect to $L>\tilde{L}_*$. If $0<L\leq \tilde{L}_*$, $\tilde{\lambda}_1^R(L)>0$ for all $R>0$,
and thus $\tilde{R}^*(L)\to +\infty$ as $L\to \tilde{L}_*+0$.

We know that the following auxiliary eigenvalue problem:
\[
\left\{
\begin{array}{ll}
- \psi'' -f'(0)\psi =\hat{\lambda} \psi, & L_1<x<L_2,\\
\psi(L_1) = \psi(L_2)=0,
\end{array}
\right.
\]
admits a principal eigenvalue $\hat{\lambda}_1$, and $\hat{\lambda}_1=0$ if $L_2-L_1=\frac{\pi}{\sqrt{f'(0)}}$.
The comparison principle implies that $\tilde{\lambda}_1^R(L)<0$ when $L\geq \frac{\pi}{\sqrt{f'(0)}}$,
which yields that $\tilde{R}^*(L)-L_1< \frac{\pi}{\sqrt{f'(0)}}$ for $L\geq \frac{\pi}{\sqrt{f'(0)}}$.

Consequently, we have
\[
\lim_{L\to \tilde{L}_*+0}\tilde{\mathcal{J}}(L)=+\infty,\ \ \ \tilde{\mathcal{J}}\Big(\frac{\pi}{\sqrt{f'(0)}}\Big)<0.
\]
Thus there exists a unique $\tilde{L}_{**}\in\Big(\tilde{L}_*,\frac{\pi}{\sqrt{f'(0)}}\Big)$ such that $\tilde{\mathcal{J}}(L_{**})=0$,
that is, $\tilde{R}^*(L_{**})=L_2$. The monotonicity of $\tilde{\mathcal{J}}(L)$ in $L$ derives all
the other assertions of the lemma.

Finally, let us give the calculation of \eqref{ee1}. It follows from the definition of $\tilde{L}_{**}$
that $\tilde{\lambda}_1^R(\tilde{L}_{**})=0$ and $\tilde{R}^*:=\tilde{R}^*(\tilde{L}_{**})=L_2$. Thus when $L=\tilde{L}_{**}$, we have
\begin{equation}\label{eigen-qr}
\left\{
 \begin{array}{ll}
 - \varphi'' -f'(0)\varphi =0, & x \in(L_1, \tilde{R}^*),\\
 - \varphi'' -g'(0)\varphi =0, & x\in(0, L_1),\\
 \varphi'(0)=\varphi(\tilde{R}^*)=0, \\
 \varphi(L_1-0) = \varphi( L_1+0),\\
 \varphi'(L_1-0) = \varphi'( L_1+0).
 \end{array}
 \right.
\end{equation}
For $x\in(0, L_1)$, since $g'(0)<0$, it follows from the second equation of \eqref{eigen-qr} that
there are two constants $\tilde{C}_1$ and $\tilde{C}_2$ such that
\[
\varphi(x)=\tilde{C}_1e^{\sqrt{-g'(0)} x}+\tilde{C}_2e^{-\sqrt{-g'(0)} x},\ \ \forall x\in (0,L_1).
\]
This, together with $\varphi'(0)=0$, yields that $\tilde{C}_1=\tilde{C}_2>0$, and so for $x\in(0,L_1)$,
\[
\varphi(x)=\tilde{C}_1\big(e^{\sqrt{-g'(0)} x}+e^{-\sqrt{-g'(0)} x}\big),
\]
\[
\varphi'(x)=\tilde{C}_1\sqrt{-g'(0)}\big(e^{\sqrt{-g'(0)} x}-e^{-\sqrt{-g'(0)} x}\big)>0.
\]
Hence we have
\begin{equation}\label{eeqq1}
\frac{\varphi'(L_1-0)}{\varphi(L_1-0)}=\sqrt{-g'(0)}\cdot\frac{e^{\sqrt{-g'(0)} L_1}
-e^{-\sqrt{-g'(0)} L_1}}{e^{\sqrt{-g'(0)} L_1}+e^{-\sqrt{-g'(0)} L_1}}.
\end{equation}
Since $\varphi'(L_1+0)=\varphi'(L_1-0)$, it follows that $\varphi'(L_1+0)>0$. As $\varphi(x)>0$ for
$x\in(0,\tilde{R}^*)$ and $\varphi(\tilde{R}^*)=0$, then there is a constant $a\in(L_1,\tilde{R}^*)$
such that $\varphi'(a)=0$.

We further claim that there exists a unique $a\in(L_1,\tilde{R}^*)$ satisfying $\varphi'(a)=0$.
Otherwise there are two constants $a_i\in (L_1,\tilde{R}^*)$ with $a_1<a_2$ satisfying $\varphi'(a_i)=0$ $(i=1,2)$.
It follows from the first equation of  \eqref{eigen-q} that
\[
-\varphi''=f'(0)\varphi, \ \ \forall x\in[a_1,a_2].
\]
Integrating the above equation on $[a_1,a_2]$ implies that $f'(0)=0$ in virtue of $\varphi'(a_i)=0$
and $\varphi>0$ in $[a_1,a_2]$, a contradiction!
Thus there is a unique $a\in(L_1,\tilde{R}^*)$ fulfilling $\varphi'(a)=0$.

Now, when $x\in(L_1,\tilde{R}^*)$, from the second equation of \eqref{eigen-qr} it follows that
there exist two constants $\tilde{C}_3$ and $\tilde{C}_4$  such that
\[
\varphi(x)=\tilde{C}_3\cos [\sqrt{f'(0)}(x-a)] +\tilde{C}_4\sin [\sqrt{f'(0)}(x-a)], \ \ \forall x\in[L_1,\tilde{R}^*].
\]
Since $\varphi'(a)=0<\varphi(x)$ for $x\in(0, \tilde{R}^*)$, then $\tilde{C}_3>0=\tilde{C}_4$, which yields that
\[
\varphi(x)=\tilde{C}_3\cos [\sqrt{f'(0)}(x-a)], \ \ \forall x\in[L_1,\tilde{R}^*].
\]
This, together with $\varphi(\tilde{R}^*)=0$, yields that
\begin{equation}\label{rr1}
\tilde{R}^*-a=\frac{\pi}{2\sqrt{f'(0)}}.
\end{equation}
Moreover, basic computation gives that
\[
\frac{\varphi'(L_1+0)}{\varphi(L_1+0)}=\sqrt{f'(0)}\tan[\sqrt{f'(0)}(a-L_1)].
\]
By virtue of \eqref{eeqq1}, it then follows that
\begin{equation}\label{rr12}
\sqrt{-\frac{g'(0)}{f'(0)}}\cdot\frac{e^{\sqrt{-g'(0)} L_1}
-e^{-\sqrt{-g'(0)} L_1}}{e^{\sqrt{-g'(0)} L_1}+e^{-\sqrt{-g'(0)} L_1}}
=\tan[\sqrt{f'(0)}(a-L_1)]>0.
\end{equation}
Clearly \eqref{ee1} follows from \eqref{rr1}, \eqref{rr12} and the definition of $\tilde{R}^*$, which ends the proof.
\end{proof}

Finally, let us give

\smallskip

\noindent
 {\bf Proof of Theorem \ref{thm:dybe1}:}\
With the help of Lemmas \ref{lem:1eigenvaluemp2} and \ref{lem:s1ei2}, one can use similar arguments as those of
Theorems \ref{thm:dybe} and \ref{thm:profile of spreading sol} (with slight modifications) to prove Theorem \ref{thm:dybe1}.
Thus all the details are omitted here. {\hfill $\Box$}
\section{Discussion}
\subsection{Analytical results}

In the current paper, we have proposed a reaction-diffusion model with strong Allee effect and free boundary, which includes
a bounded protection zone, in order to save the endangered species. We employ the contraction mapping theorem, the comparison
principle and $L^p$ estimates together with the Sobolev embedding theorem to establish the global existence and uniqueness of
solutions of problems \eqref{p} and \eqref{q}. Then we find that there exist three critical values $\overline{L}_*<\overline{L}_{**}$
and $\overline{L}^*$, which play key roles in the dynamics of the solutions and divide the protection zone into three cases: small
protection zone case ($0< L\leq \overline{L}_*$); medium-sized protection zone case ($\overline{L}_* < L < \max\{\overline{L}^*, \overline{L}_{**}\}$)
and large protection zone case ($L>\max\{\overline{L}^*, \overline{L}_{**}\}$), see Theorem \ref{thm:dybe} ($\overline{L}_*=L_*$,
$\overline{L}^*=L^*$,  $\overline{L}_{**}=L_{**}$) and Theorem \ref{thm:dybe1} ($\overline{L}_*=\tilde{L}_*$, $\overline{L}^*=\tilde{L}^*$,
$\overline{L}_{**}=\tilde{L}_{**}$). In the small protection zone case ($0<L\leq \overline{L}_*$), there is a vanishing-transition-spreading
trichotomy result. In the medium-sized protection zone case ($\overline{L}_* <L< \max\{\overline{L}^*, \overline{L}_{**}\}$),
if $\overline{L}^*< \overline{L}_{**}$, there is a vanishing-transition-spreading trichotomy result when $\overline{L}_*<L<\overline{L}^*$;
while there is a vanishing-spreading dichotomy result when $\overline{L}^*<L<\overline{L}_{**}$, and if $\overline{L}^*>\overline{L}_{**}$,
there is a vanishing-transition-spreading trichotomy result when $\overline{L}_*<L<\overline{L}_{**}$; while there is a transition-spreading
dichotomy result when $\overline{L}_{**}<L<\overline{L}^{*}$. In the large protection zone case ($L>\max\{\overline{L}^*,\overline{L}_{**}\}$),
only spreading happens. As a consequence, our results reveal that the species will survive in the entire space if the length of the protection
zone is larger than $\overline{L}_{**}$, which is called the effective length of the protection zone.  Moreover, (iii) of Lemma \ref{lem:s1ei2} tells us that $\tilde{L}_{**}>L_{**}$. Combining
this with Theorems \ref{thm:dybe} and \ref{thm:dybe1}, we conclude that the connected protection zone is better for species survival than a
separated one. Furthermore, when spreading happens, a precise estimate for the spreading speed and the asymptotic profiles are given in Theorems
\ref{thm:profile of spreading sol} and \ref{thm:dybe1}.

\subsection{Discussion and conclusions}

Our Theorem \ref{thm:dybe} suggest that
the asymptotic profiles of the evolution of the endangered species behave differently from that of the corresponding problem \eqref{pp1}
without free boundary condition and that of the problem \eqref{20201031-2} without the  protection zone. In the following, we would like to make a
comparison and reveal the effect of the protection zone and free boundary on the spatial distribution of the endangered species.

\smallskip

\subsubsection{The effect of the protection zone on the asymptotic behavior of solutions of \eqref{p}.}
Comparing to the results on the solution of problem \eqref{20201031-2} (in \cite{DuLou} Theorem 1.3), we see that the spatial distribution of the solution of \eqref{p} has the following four cases: either a
spreading-transition-vanishing trichotomy result, or a spreading-vanishing dichotomy result, or a spreading-transition dichotomy result, or a spreading
phenomenon holds (see Theorem \ref{thm:dybe} for more details). Moreover, the long-time dynamical behaviors of
the solution of \eqref{p} do not vanish  provided that the protection zone is
suitably large ($L>{L}_{**}$). In a biological sense, the endangered species will survive in the whole habitat
regardless of the initial density of species. This implies that the protection zone can enhance persistence of the endangered species.

\smallskip

\subsubsection{The effect of the free boundary condition on the asymptotic behavior of solutions of \eqref{p}}

In the point of saving the endangered species, problems \eqref{pp1} and \eqref{p} have introduced a protection zone ($0<x<L$). As for problem \eqref{p},
there is another critical value $L_{**}(>L_*)$ which is determined by the free boundary condition, besides the same two critical values $L_*$ and $L^*$
in problem \eqref{pp1}. The third critical value $L_{**}$ divides the long-time dynamical behavior of the solutions for problem \eqref{p} into two types in
the medium-sized protection zone. The first one is either a spreading-transition-vanishing trichotomy result or a spreading-vanishing dichotomy result;
and the second one is either a spreading-transition-vanishing trichotomy result or a spreading-transition dichotomy result, please see Theorem \ref{thm:dybe}.
But there is only a spreading-transition  dichotomy result for problem \eqref{pp1} in the medium-sized protection zone in \cite{DPS}.

Now, we consider the effective length of the protection zone of problems \eqref{p} and \eqref{pp1}. It follows from Theorem 1.1 in \cite{DPS} that the asymptotic
profiles of the solution for problem \eqref{pp1} do not vanish for any initial value as long as $L>L_*$.  But for the problem \eqref{p}, our results in Theorem
\ref{thm:dybe} suggest that either a spreading-transition dichotomy result or a spreading phenomenon happens for problem \eqref{p} when $L>L_{**}(>L_*)$. In other
words, the effective length of the protection zone of problem \eqref{p} is longer than that of problem \eqref{pp1}. As a result, the presence of the free boundary
condition increases the effective length of the protection zone, which yields that the free boundary condition makes the endangered species harder to survival
than that without one.

\smallskip

From what has been discussed above, we may safely draw the conclusion that both of the free boundary
condition and the protection zone play a significant role on the survival of the endangered species, and make the dynamics of solutions much more complicated and meaningful than that in \cite{DPS,DuLou}.

\end{document}